\documentclass[a4paper, reqno]{amsart}
\usepackage{amssymb}
\usepackage{amsmath}
\usepackage{amsthm}
\usepackage[utf8]{inputenc}
\usepackage[T1]{fontenc}
\usepackage{lmodern}
\usepackage{ae}
\usepackage{graphicx}
\usepackage{enumerate}
\usepackage{fancybox}
\usepackage{setspace}
\usepackage[all,xdvi]{xy}
\usepackage{ifthen}
\usepackage{everypage}
\usepackage{tikz}
\usetikzlibrary{shapes,matrix,arrows,calc}

\title[Cyclotomic associators and tangles invariants in the torus]{Cyclotomic associators and finite type invariants for tangles in the solid torus}
\author{Adrien Brochier}
\address{Université de Genève}
\email{adrien.brochier@unige.ch}

\tikzstyle cross=[preaction={draw=white, -, line width=4pt}, thick]
\tikzstyle normal=[thick]
\tikzstyle chord=[densely dotted, thick]
\tikzstyle zero=[ultra thick, gray]
\tikzstyle zerocross=[preaction={draw=white, -, line width=4pt}, ultra thick, gray]
\tikzstyle point=[draw,circle,inner sep=1,fill=black]
\tikzstyle petitpoint=[draw,circle,inner sep=0.3,fill=black]
\def\taosize {0.65}
\newcommand{\negative}[3][-]{\draw[normal,#1] (#2+1,-#3).. controls (#2+1,-#3-0.3) and (#2,-#3-0.7)..(#2,-#3-1); \draw[cross,#1] (#2,-#3).. controls (#2,-#3-0.3) and (#2+1,-#3-0.7)..(#2+1,-#3-1);}
\newcommand{\positive}[3][-]{ \draw[normal,#1] (#2,-#3).. controls (#2,-#3-0.3) and (#2+1,-#3-0.7)..(#2+1,-#3-1);\draw[cross,#1] (#2+1,-#3).. controls (#2+1,-#3-0.3) and (#2,-#3-0.7)..(#2,-#3-1);}
\newcommand{\sing}[3][-]{ \draw[normal,#1] (#2,-#3).. controls (#2,-#3-0.3) and (#2+1,-#3-0.7)..(#2+1,-#3-1);\draw[normal,#1] (#2+1,-#3).. controls (#2+1,-#3-0.3) and (#2,-#3-0.7)..(#2,-#3-1);\node[point] at (#2+0.5,-#3-0.5) {}; }
\newcommand{\perm}[3][-]{ \draw[normal,#1] (#2,-#3).. controls (#2,-#3-0.3) and (#2+1,-#3-0.7)..(#2+1,-#3-1);\draw[normal,#1] (#2+1,-#3).. controls (#2+1,-#3-0.3) and (#2,-#3-0.7)..(#2,-#3-1); }
\newcommand{\twist}[3][-]{
	\begin{scope}[xshift=#2cm,yshift=-#3cm,yscale=0.5,xscale=0.8]
		\draw[normal] (0,0)..controls (0,-1) and (0.5,-1).. (0.5,-0.5);
		\draw[cross,#1] (0.5,-0.5)..controls (0.5,0) and (0,0).. (0,-1);
	\end{scope}
}
\newcommand{\twistinv}[3][-]{
	\begin{scope}[xshift=#2cm,yshift=-#3cm,yscale=0.5,xscale=0.8]
		\draw[normal,#1] (0.5,-0.5)..controls (0.5,0) and (0,0).. (0,-1);
		\draw[cross] (0,0)..controls (0,-1) and (0.5,-1).. (0.5,-0.5);
	\end{scope}
}
\newcommand{\tainv}[2][-]{ 
\begin{scope}[xscale=\taosize]
\draw[normal] (0,-#2).. controls (0,-#2-0.2) and (0-1.2,-#2-0.2)..(0-1.2,-#2-0.5);
\draw[zerocross] (-1,-#2)--(-1,-#2-0.5);
\draw[zero] (-1,-#2-0.5)--(-1,-#2-1);
\draw[cross,#1] (-1.2,-#2-0.5).. controls (-1.2,-#2-0.7) and (0,-#2-0.7)..(0,-#2-1);
\end{scope}
}
\newcommand{\taoo}[2][-]{ 
\begin{scope}[xscale=\taosize]
\draw[zero] (-1,-#2)--(-1,-#2-0.5);
\draw[cross] (0,-#2).. controls (0,-#2-0.2) and (-1.2,-#2-0.2)..(-1.2,-#2-0.5);
\draw[normal,#1] (-1.2,-#2-0.5).. controls (-1.2,-#2-0.7) and (0,-#2-0.7)..(0,-#2-1);
\draw[zerocross] (-1,-#2-0.5)--(-1,-#2-1);
\end{scope}
}
\newcommand{\taotet}[2][-]{
	\begin{scope}[yshift=-#2cm,yscale=0.666666666]
	\taoo{0}
	\fixed{0.5}
	\twist[#1]{0}{1}
	\end{scope}
}
\newcommand{\taotetinv}[2][-]{
	\begin{scope}[yshift=-#2cm,yscale=0.666666666]
	\tainv{0}
	\fixed{0.5}
	\twistinv[#1]{0}{1}
	\end{scope}
}

\newcommand{\staoo}[2][-]{ 
\begin{scope}[xscale=\taosize]
\draw[zero] (-1,-#2)--(-1,-#2-0.5);
\draw[cross] (0,-#2).. controls (0,-#2-0.2) and (-1.2,-#2-0.2)..(-1.2,-#2-0.5);
\draw[normal,#1] (-1.2,-#2-0.5).. controls (-1.2,-#2-0.7) and (0,-#2-0.7)..(0,-#2-1);
\draw[zero] (-1,-#2-0.5)--(-1,-#2-1);
\node[point] at (-1,-#2-0.63) {} ;
\end{scope}
}

\newcommand{\up}[3][-]{ 
\draw[normal,#1] (#2,-#3).. controls (#2,-#3-0.3) and (#2+0.2,-#3-0.5).. (#2+0.5,-#3-0.5);
\draw[normal] (#2+0.5,-#3-0.5).. controls (#2+0.8,-#3-0.5) and (#2+1,-#3-0.3).. (#2+1,-#3);}
\newcommand{\ap}[3][-]{
\draw[normal,#1] (#2,-#3-1).. controls (#2,-#3-0.7) and (#2+0.2,-#3-0.5).. (#2+0.5,-#3-0.5);
\draw[normal] (#2+0.5,-#3-0.5).. controls (#2+0.8,-#3-0.5) and (#2+1,-#3-0.7).. (#2+1,-#3-1);}

\newcommand{\straight}[3][-]{\draw[normal,#1] (#2,-#3) -- (#2,-#3-1);}
\newcommand{\zn}[4][-]{\draw[normal,#1] (#2,-#3) -- (#2,-#3-1);\node[point,label=left:$#4$] at (#2,-#3-0.5) {};}
\newcommand{\fixed}[1]{\draw[zero] (-\taosize,-#1) -- (-\taosize,-#1-1);}
\newcommand{\fixch}[1]{\draw[zero] (0,-#1) -- (0,-#1-1);}

\newcommand{\sdot}[2]{\draw[dotted] (#1,-#2) -- (#1,-#2-1);}
\newcommand{\conf}[2]{ 
\foreach \angle in {0,72,144,216,288} {\node[petitpoint] at (\angle:#1) {}; }
\node[petitpoint,label=below:$0$] at (0,0) {};
\node[label=above:$\zeta_N$] at (72:#1) {};
\node[label=below:$1$] at (1,0) {};
\draw[dotted] (-#2,-#2) rectangle (#2,#2) ;
}
\newcommand{\hori}[4][-]{\draw[normal,#1] (#2,-#3)--(#2,-#3-1);\draw[normal,#1] (#2+#4,-#3)--(#2+#4,-#3-1);\draw[chord] (#2,-#3-0.5)--(#2+#4,-#3-0.5);}
\newcommand{\labhori}[6][-]{\draw[normal,#1] (#2,-#3)--(#2,-#3-1);\draw[normal,#1] (#2+#4,-#3)--(#2+#4,-#3-1);\draw[chord] (#2,-#3-0.5)--(#2+#4,-#3-0.5);
\node[point,label=left:$#5$] at (#2,-#3) {};\node[point,label=left:$#6$] at (#2,-#3-1) {};
}
\newcommand{\tzero}[3][-]{\draw[zero] (0,-#2)--(0,-#2-1);\draw[normal,#1] (#3,-#2)--(#3,-#2-1);\draw[chord] (0,-#2-0.5)--(#3,-#2-0.5);}
\newcommand{\tik}[1]{\begin{tikzpicture}[baseline=(current bounding box.center)] #1 \end{tikzpicture} }
 \usepackage[pdftex,
     bookmarks         = true,
     bookmarksnumbered = true,
     pdfpagemode       = None,
     pdfstartview      = FitH,
     pdfpagelayout     = SinglePage,
     urlcolor          = blue,
     pdfborder         = {0 0 0},
 	pdftoolbar = true
     ]{hyperref}                

\numberwithin{equation}{section}

\newcommand{\RR}{\mathbb{R}}
\newcommand{\C}{\mathbb{C}}

\newcommand{\Z}{\mathbb{Z}}
\newcommand{\Q}{\mathbb{Q}}

\newcommand{\cd}{\operatorname{\bf CD}}
\newcommand{\bcd}{\operatorname{\bf BCD}}

\newcommand{\cyc}{\operatorname{Cyclo}}

\newcommand{\mf}[1]{\mathfrak{#1}}

\newcommand{\id}{\operatorname{id}}

\newcommand{\im}{\operatorname{Im}}

\newcommand{\R}{\mathcal{R}}

\newcommand{\h}{\hslash}

\newcommand{\ot}{\otimes}

\newcommand{\T}{\mathcal T}
\newcommand{\ST}{\mathcal {ST}}
\newcommand{\I}{\mathcal I}

\newcommand{\mC}{\mathcal {C}}
\newcommand{\mM}{\mathcal {M}}
\newcommand{\kk}{\mathbb K}
\newcommand{\tang}{\operatorname{\bf Tang}}
\newcommand{\btang}{\operatorname{\bf BTang}}
\renewcommand{\hom}{\operatorname{Hom}}
\newcommand{\eno}{\operatorname{End}}

\newcommand{\aut}{\operatorname{Aut}}

\newcommand{\modu}{\operatorname{-mod}}

\newcommand{\gr}{\operatorname{gr}}
\newcommand{\src}{\operatorname{source}}
\newcommand{\tgt}{\operatorname{target}}
\newcommand{\coev}{\operatorname{coev}}
\newcommand{\ev}{\operatorname{ev}}
\newtheorem{thm}{Theorem}[section]
\newtheorem{cor}[thm]{Corollary}
\newtheorem{prop}[thm]{Proposition}
\newtheorem{conj}[thm]{Conjecture}

\let\oldcite\cite
\renewcommand{\cite}[2][noOption]{\ifthenelse{\equal{#1}{noOption}}{{\scriptsize \oldcite{#2}}}{{\scriptsize \oldcite[#1]{#2}}}}
\newtheorem{defi}[thm]{Definition}

\theoremstyle{remark} 
\newtheorem{rmk}[thm]{Remark}

\begin{document}
\begin{abstract}
The universal Vassiliev-Kontsevich invariant is a functor from the category of tangles to a certain graded category of chord diagrams, compatible with the Vassiliev filtration and whose associated graded is an isomorphism. The Vassiliev filtration has a natural extension to tangles in any thickened surface $M\times I$ but the corresponding category of diagrams lacks some finiteness properties which are essential to the above construction. We suggest to overcome this obstruction by studying families of Vassiliev invariants which, roughly, are associated to finite coverings of $M$. In the case $M=\C^*$, it leads for each positive integer $N$ to a filtration on the space of tangles in $\C^* \times I$ (or "B-tangles"). We first prove an extension of the Shum--Reshetikhin--Turaev theorem in the framework of braided module category leading to B-tangles invariants. We introduce a category of "$N$-chord diagrams", and use a cyclotomic generalization of Drinfeld associators, introduced by Enriquez, to put a braided module category structure on it. We show that the corresponding functor from the category of B-tangles is a universal invariant with respect to the $N$ filtration. We show that Vassiliev invariants in the usual sense are well approximated by $N$ finite type invariants. We show that specializations of the universal invariant can be constructed from modules over a metrizable Lie algebra equipped with a finite order automorphism preserving the metric. In the case the latter is a "Cartan" automorphism, we use a previous work of the author to compute these invariants explicitly using quantum groups. Restricted to links, this construction provides polynomial invariants.
\end{abstract}

\maketitle
\setcounter{section}{-1}
\section{Introduction}
Finite type invariants are those numerical knot invariants whose "$(n+1)$th derivative", in some specific sense, vanishes for some $n$. Since their discovery around 1989 by Vassiliev, it became clear that there is a deep relation between finite type invariants, Lie theory and quantum topology. Indeed, the collection of all finite type invariants dominates all quantum invariants, including the various knot polynomials. Kontsevich~\cite{Kontsevich1993} then gave an essentially complete description of the space of finite type invariants, by constructing a "universal" invariant taking its values in some space of Feynman diagrams. The Kontsevich integral not only provide a powerful knot invariant, but shed some light on the topological background of deformation-quantization theory.

Let $\kk$ be a field of characteristic 0 and $V$ be the $\kk$-linear span of isotopy classes of knots in $\RR^3$, then $\kk$-valued knots invariants may be identified with linear map $V\rightarrow \kk$. A singular knot can be identified with a formal linear combination of non-singular knots by repetitive use of the Vassiliev skein relation\footnote{Here as usual we draw only the part of the knots which is involved in the relation, and assume that they are identical outside the picture.}~\cite{Vassiliev1990}:
\[
\tik{\sing[->]{0}{0}}=\tik{\positive[->]{0}{0}}-\tik{\negative[->]{0}{0}}
\]
This define a filtration on $V$ by letting $\I_n$ to be the linear span of knots having a least $n$ singularities. A \emph{finite type} invariant of type (or degree) at most $n$ is a knot invariant whose extension to $V$ vanishes on $\I_{n+1}$ for some $n$, i.e. an element of $\hom(V/\I_{n+1},\kk)$.

It follows that an invariant of type $n$ applied to a singular knot with $n$ singularities is blind to its topology, meaning that its value does not change if one replace an overcrossing by an undercrossing and \emph{vice versa}. Hence, it is sensitive only to the \emph{combinatorial} information corresponding to the position of singularities. Therefore, any finite type invariant induces an invariant of a combinatorial object (a chord diagram) encoding this information. One then shows that every finite type invariant descend to a linear map on the space $A$ of chord diagrams modulo the so-called 4T relation (see Section~\ref{sec:funda}). Since the latter is homogeneous, $A$ is graded and hence inherits a natural filtration. The fundamental theorem of Vassiliev invariants states that the space of finite type invariants is isomorphic, as a filtered vector space, to the graded dual of $A$. We refer the reader to Bar-Natan's survey~\cite{Bar-Natan1995} for an introduction on finite type invariants, and to Chmutov--Duzhin--Mostovoy's book~\cite{Chmutov2012} for a fairly complete account of the subject.

The proof of this result relies on the construction of a \emph{universal} finite type invariant: this is a filtration preserving linear map $Z$ from $V$ to the degree completion of $A$, whose associated graded is an isomorphism. Hence, if $f$ is a finite type invariant and $\tilde f$ the linear form it induces on $A$ (its "symbol") then
\[
	f=\tilde f \circ Z +(\text{ lower degree terms }).
\]

This was first proved for knots by Kontsevich~\cite{Kontsevich1993} by integrating a formal analog of the KZ connection. As hinted by Kontsevich himself, it was then observed that such an invariant can be obtained through a kind of Markov trace from Drinfeld's universal representation of the braid group, which also comes from the KZ connection~\cite{Drinfeld1990a}. This representation depends algebraically on a formal power series (the KZ associator) in two non-commuting variables satisfying a set of algebraic equations. It turns out that any formal power series satisfying the same set of equations leads to a representation of the braid group having similar properties. Such elements are called Drinfeld associators.

Hence, a combinatorial formula for a universal finite type invariant involving Drinfeld associators was then given by Le--Murakami~\cite{Le1995} and Piunikhin~\cite{Piunikhin1995}, allowing the generalization of Kontsevich's theorem to tangles. Finally, it was put in the categorical language by Cartier~\cite{Cartier1993}, Bar-Natan~\cite{Bar-Natan1997} and Kassel--Turaev~\cite{Kassel1998}. 

The proof goes as follows: both tangles and chord diagrams can be thought of as morphisms in certain categories $\tang$ and $\cd$ having the same objects. It follows from~\cite{Drinfeld1990a} that each choice of a Drinfeld associator leads to a ribbon structure on the degree completion $\widehat{\cd}$ of $\cd$. By Shum's theorem~\cite{Shum1994}, the tangle category is universal among ribbon categories, hence one gets a functor
\[
F:\tang \longrightarrow \widehat{\cd}
\]
which is shown to be a universal finite type invariant. Working in a categorical setting already makes it clear that this invariant is compatible with the composition of tangles whenever it is defined. Drinfeld prove in \emph{loc. cit.} that associators with rational coefficients exist and can be constructed recursively. Hence such a functor exists for $\kk=\Q$ and can be combinatorially computed up to any degree in a finite amount of time. 

A very important feature of the space of invariants of some given degree (or equivalently of the graded pieces of the space of diagrams) is that it is finite dimensional. In particular, the space of diagrams associated to pure braids is a finitely generated algebra, which is the basic ingredient in the definition of the KZ equation and of Drinfeld associators.

While the Vassiliev skein relation makes perfectly sense for knotted objects in a thickened surface $M\times I$, the above mentioned finiteness properties fails if $\pi_1(M)$ is infinite. This can be roughly explained as follows: it is still true that a degree $n$ finite type invariant applied to a knot $K$ in $M\times I$ with $n$ singularities is blind to the topology of $K$. But it is, in general, not blind to the topology of the surface itself, in that in can see the element of $\pi_1(M)$ induced by $K$. Hence, it induces a function on a space of diagrams labelled by $\pi_1(M)$ whose graded pieces are infinite dimensional.

This is a strong obstruction to the construction of a universal finite type invariant enjoying some desirable features. It implies for example that one cannot use Chen iterated integrals which are the basic ingredient of Kohno's construction of a universal finite type invariant for pure braids~\cite{Kohno1983}. Indeed, it is known~\cite{Bellingeri2004157} that if the genus of $M$ is strictly positive, then a universal \emph{multiplicative} invariant for the pure braid groups of $M$ cannot exist. A universal invariant for links thus cannot be obtained through a Markov trace.

This motivates the study of sub-families of finite type invariant which, when applied to a knot (or tangle) with enough singularity, distinguishes its representative in $\pi_1(M)$ only modulo a given cofinite subgroup $H$ of $\pi_1(M)$. We give a general definition in the case of the pure braid group $P_n(M)$ of $M$ in Section~\ref{sec:motiv}, explain how this is related to finite covering of $M$ and show that in many cases the collection of these invariants is as strong as the collection of finite type invariants. That $H$ is cofinite implies that the associated graded of $P_n(M)$ with respect to this filtration is finitely generated. This definition also makes sense if $H$ is not cofinite, and the case $H=\{1\}$ gives back the Vassiliev filtration. On the other hand, the filtration associated with $H=\pi_1(M)$ coincides with the filtration by the powers of the augmentation ideal of $\kk[P_n(M)]$, which was extended to tangles by Habiro under the name "A-filtration"~\cite{Habiro2000}. Our construction can thus be thought of as an interpolation between theses two filtrations, which coincide in the classical case.

We then focus on the case $M=\C^*$, i.e. on tangles in the solid torus (which we call B-tangles). Note that a universal finite type invariant for B-knot~\cite{Goryunov1997} and B-links~\cite{Andersen1998,Lieberum2004} was already constructed, but they are \emph{a priori} not functorial nor extensible to tangles. Also, that they take values in a space of infinite type makes them difficult to handle from an algorithmic point of view, contrasting with the classical case.

Recall that the skeleton of a tangle is its underlying abstract oriented 1-manifold. We introduce for each integer $N\geq 1$ a $N$-filtration $\{\I_{n,N}(S)\}_{n\geq 1}$ on the space of B-tangles on a given skeleton $S$, which has a natural geometric interpretation by extending the notion of singular tangles. It leads to a notion of $N$-finite type invariants (section~\ref{sec:def}). As a first example, we show that, upon some specialization, the invariants for B-knots constructed by Lambropoulou in~\cite{Lambropoulou1999} from traces on cyclotomic Hecke algebras are of $N$-finite type (section~\ref{sec:hecke}).

We then introduce (section~\ref{sec:ndiag}) a space of $N$-diagrams and show that any $N$-finite type invariant descend to a function on the quotient $B_N(S)$ of the space of diagrams by a set of explicit relations. 

The main result of this paper can be stated as follows:
\begin{thm}
For all $N\geq 1$, for any skeleton $S$, there exists a well-behaved  universal $N$-finite type invariant, that is a combinatorial, filtration preserving, functorial linear map defined over $\Q$ from the space of linear combinations of B-tangles with skeleton $S$ to $B_N(S)$, whose associated graded is an isomorphism.
\end{thm}
In particular, the space of $N$-finite type invariants is isomorphic as a filtered vector space to the graded dual of $B_N(S)$.

We recall in Section~\ref{sec:kv} the categorical construction of the Vassiliev--Kontsevich invariant. We prove in Section~\ref{sec:shum} a generalisation of Shum's theorem in the framework of braided module category. It provides a systematic way of producing B-tangles invariants. Key to the second part of the proof is a cyclotomic analog of the notion of associator developed by Enriquez~\cite{Enriquez2008}: this is a pair $(\Psi,\Phi)$ where $\Phi$ is an ordinary associator, and $\Psi$ is a formal power series in $N+1$ non-commuting variables satisfying a set of equations. We give our main construction in Section~\ref{sec:main}: we introduce a notion of $N$ infinitesimal braided module category and show that the category $\bcd_N$ of $N$ diagrams is universal among such structures. We prove that the choice of a cyclotomic associator $(\Psi,\Phi)$ turns the degree completion $\widehat{\bcd}_N$ of $\bcd_N$ into a braided module category over the ribbon category structure on $\widehat{\cd}$ coming from $\Phi$. Finally, we show that the functor
\[
G:\btang\longrightarrow \widehat{\bcd}_N
\]
is a universal $N$-finite type invariant.

We then show in Section~\ref{sec:approx} that usual Vassiliev invariants are "well approximated" by the collection of all $N$-finite type invariants. More precisely, we construct a functor from $\btang$ to the "profinite completion" of the categories $\bcd_N$ for all $N$, whose associated graded is faithful.

Finally, in Section~\ref{sec:quant}, we construct specializations of this invariant (a.k.a. weight systems). It is well known that representations of the category of ordinary diagrams can be constructed out of a metrizable Lie algebra, that is, a Lie algebra $\mf g$ together with a non-degenerate $t \in S^2(\mf g)^{\mf g}$. In very much the same way, equipping $\mf g$ with a $N$th order automorphism $\sigma$ preserving the metric allows to construct realizations of $\bcd_N$ for each $\mf g$-module $V$ and each $\mf h$-module $M$ where $\mf h=\mf g^{\sigma}$. Using author's previous results~\cite{Brochier2012}, we show that these specializations can be explicitly computed using quantum groups, in the case where $\mf g$ is simple and $\sigma$ is a Cartan automorphism. Restricted to links and if $M$ is a one dimensional $\mf h$-module given by an integral weight, these invariants are, in fact, Laurent polynomials in one variable, generalizing the Jones polynomial and other polynomial invariants in a very natural way.

\medskip
\noindent
{\bf Acknowledgments.} I would like to thank Benjamin Enriquez, Philippe Humbert and Gw\'ena\"el Massuyeau for useful discussions and comments. I also would like to thank the referee for many helpful remarks which improved the quality of this paper.
\section{Reminder and definitions}\label{sec:def}
\subsection{Motivation: $H$-finite type invariants for pure braids}\label{sec:motiv}
Let $M$ be an orientable surface. The pure braid group $P_n(M)$ of $M$ is the fundamental group of the configuration space
\[
X_n(M)=\{(z_1,\dots,z_n) \in M^n,i\neq j \Rightarrow z_i \neq z_j \}
\]
There is a canonical surjective group morphism
\[
P_n(M)\longrightarrow \pi_1(M)^n
\]
which maps any pure braid to the $n$-uplet formed by the homotopy classes of its $n$ strands. It induces an algebra map
\[
	\kk[P_n(M)]\longrightarrow \kk[\pi_1(M)^n].
\]
Let $I_n$ be its kernel, then the Vassiliev filtration for $\kk[P_n(M)]$ coincides with the filtration induced by the powers of $I_n$.

If $p:\tilde M\rightarrow M$ is a regular finite covering of $M$ associated to a cofinite subgroup $H$ of $\pi_1(M)$, one can define the orbit configuration space
\[
X_{n,H}(M)=\{ (z_1,\dots, z_n) \in \tilde M, i \neq j \Rightarrow p(z_i)\neq p(z_j) \}
\]
Let $P_{n,H}(M)$ be the fundamental group of $X_{n,H}$ and $G_H=\pi_1(M)/H$, then $p$ induces a short exact sequence
\[
1\rightarrow P_{n,H}(M) \rightarrow P_n(M) \rightarrow G_H^n \rightarrow 1
\]
where the last map is the composition
\[
P_n(M) \rightarrow \pi_1(M)^n \rightarrow G_H^n
\]

Through the natural identification $X_n(M)=X_{n,H}(M)/G_H^n$, one can represent elements of $P_n(M)$ as an union of paths in $\tilde M$ whose start point and end point are not necessarily equal but related by the action of $G_H$. 
\begin{defi}
Let $\I_{n,H}(M)$ be the kernel of the corresponding algebra map
\[
\kk[P_n(M)]\longrightarrow \kk[G_H^n]
\]

An invariant of $H$-finite type is a linear map
\[
f:\kk[P_n(M)]\rightarrow \kk
\]
such that $f(\I_{n,H}(M)^k)=0$ for some $k$.
\end{defi}
\begin{rmk}
Although it is not cofinite, it makes sense to take $H=\{1\}$ in the above definition. In that case, one obviously recover the usual notion of finite type invariant. On the other hand, if $H=\pi_1(M)$, $\I_{n,H}(M)$ is the augmentation ideal of $\kk[P_n(M)]$.
\end{rmk}

This definition is motivated by the following easy fact:
\begin{prop}
Assume that $\pi_1(M)$ is residually finite. Then for all $n\geq 0$,
\[
\bigcap_{H\text{ cofinite subgroup of } \pi_1(M)} \I_{n,H}=\I_n
\]
\end{prop}
It means that if $\pi_1(M)$ is residually finite, then the collection of all $H$-finite type invariants for all cofinite subgroups $H$ is as strong as the collection of finite type invariants in the usual sense.
Since $P_{n,H}(M)$ is a finitely generated fundamental group of a nice topological space, contrasting with the above-mentioned negative result~\cite{Bellingeri2004157}, we can make the following conjecture:
\begin{conj}
Let $I_{P_{n,H}(M)}$ be the augmentation ideal of $\kk[P_{n,H}(M)]$ and
\[
	\gr \kk[P_{n,H}(M)]=\prod_{k\geq 0} I_{P_{n,H}(M)}^{k}/I_{P_{n,H}(M)}^{k+1}.
\]
There exists a filtration preserving algebra morphism
\[
\kk[P_n(M)] \rightarrow \gr \kk[P_{n,H}(M)]\rtimes G_H^n
\]
whose associated graded is the identity.
\end{conj}
Indeed, this conjecture is true in the following cases:
\begin{itemize}
\item If $M$ is a torus, $\kk=\Q$ and $H=\Z^2$ (\cite{Calaque2009,Humbert2012});
\item If $M$ is a compact surface of genus $g$, $\kk=\C$ and $H=\pi_1(M)$ (\cite{Bezrukavnikov1994,Enriquez2005b});
\item If $M=\C^*$, $\kk=\Q$ and $H=N\Z$ (\cite{Enriquez2008});
\end{itemize}

Assume now that $M=\C^*$. In that case, the point $z=0$ can be thought of as an additional fixed strand. In particular, $P_n(\C^*)$ is isomorphic to $P_{n+1}$, the usual pure braid group on $n+1$ strands. Recall that the braid group $B_n^1$ of $\C^*$ is the braid group of Coxeter type B.
\begin{prop}[\cite{Brieskorn1971}]
 The group $B_n^1$ admits the following presentation:
\begin{align*}
B_n^1=\langle \tau, \sigma_1,\dots, \sigma_{n-1}\ |\ & \tau \sigma_1 \tau \sigma_1 = \sigma_1 \tau \sigma_1 \tau \\ 
& \tau \sigma_i=\sigma_i\tau \text{ if } i > 1 \\
 & \sigma_i \sigma_{i+1} \sigma_i = \sigma_{i+1} \sigma_{i} \sigma_{i+1}\ \forall  i \in \{1,\dots,n-2 \}\\
 & \sigma_i \sigma_j = \sigma_j \sigma_i \text{ if } |i-j| \geq 2 \rangle
\end{align*}
\end{prop}

Let $\I_{n,N}$ be the kernel of the algebra morphism
\[
\kk[P_{n+1}]\longrightarrow \kk[(\Z/N\Z)^n]
\]
and $\tilde \I_{n,N}$ be the ideal in $\kk[B_n^1]$ generated by
\begin{align*}
\tau^N&-1&\sigma-\sigma^{-1}
\end{align*}
It is easily checked that 
\[
\I_{n,N}=\tilde \I_{n,N}\cap \kk[P_{n+1}]
\]

This allows to give a description of the above filtration in terms of singular braid, by allowing singularities with this additional strand as well. Moreover, working instead with paths in the finite covering of $\C^*$ associated to the map $z\mapsto z^N$ allows to draw singular analogs of each $\tau^k$, $k=0\dots N-1$ and their corresponding resolution. For example, for $k=1$:
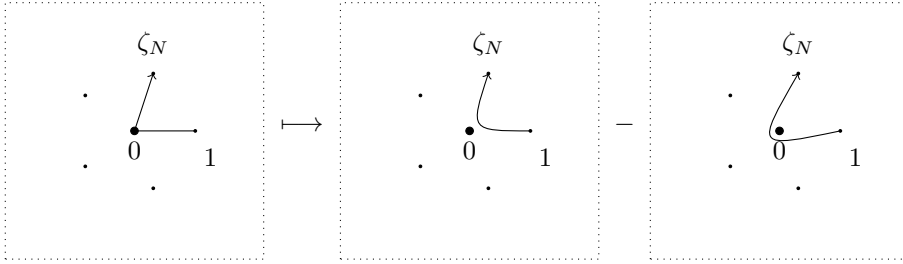
\[
\tik{\conf{0.8}{1.7} \draw[->] (0.8,0)--(0,0) -- (72:0.8); \node[point] at (0,0) {} ;  }\ \longmapsto \ \tik{\conf{0.8}{1.7} \draw[->] (0.8,0)..controls (0,0) .. (72:0.8); \node[point] at (0,0) {} ;  }\ - \ \tik{\conf{0.8}{1.7} \draw[->] (0.8,0).. controls (-144:0.4) .. (72:0.8); \node[point] at (0,0) {} ;  }
\]

It shows that, like in the usual Vassiliev theory, elements of $\tilde\I_{n,N}$ are obtained by taking the difference between the two possible small perturbations of a singular path.
\subsection{Tangles, chord diagrams and the fundamental theorem}\label{sec:funda}
Let us begin with some standard definition:
\begin{defi}
A \emph{skeleton} $S$ is the data consisting of:
\begin{itemize}
\item a one dimensional, compact, oriented manifold $X$ with (maybe empty) boundary;
\item a partition of $\partial X$ into two sets $\src(S)$ and $\tgt(S)$ together with a choice of a total ordering on these sets.
\end{itemize}
The sequence $\src(S)$ (resp.$\tgt(S)$) is represented by a word $s_1\dots s_k$ (resp. $t_1\dots t_l$) in $\{+,-\}$ where $s_i=+$ (resp. $t_i=-$) if the corresponding interval is oriented towards the corresponding boundary component and $-$ (resp. $+$) otherwise.
\end{defi}
\begin{defi}
Two skeleton $S_1,S_2$ are \emph{composable} if $\tgt(S_1)=\src(S_2)$. In that case the composition of the two is obtained by stacking $S_1$ over $S_2$.
\end{defi}
\begin{defi}
	An oriented \emph{tangle} $T$ with underlying skeleton $S$ is a smooth embedding $\iota$ of (the underlying manifold of) $S$ into $\C \times [0,1]$ mapping $s_i$ to $(i,0)$ and $t_i$ to $(i,1)$. A framing on $T$ is an extension of $\iota$ to an embedding of the "ribbon" $X\times [0,1]$, or equivalently a choice of a normal unit vector field on $\iota(X)$. We assume that $\iota(X)$ meet $\C\times \{0\}$ and $\C\times \{1\}$ orthogonally, and that the framing at each endpoint is given by $(-\sqrt{-1},0)$. Framed oriented tangles are considered up to isotopies preserving the boundary.
\end{defi}
\begin{defi}
A singular tangle with skeleton $S$ is a smooth map $X \rightarrow \C \times [0,1]$ with finitely many transverse double points, considered up to isotopy fixing singular points. Orientation and framing are defined similarly to the case of tangles.
\end{defi}

Composition of skeleta induces a composition at the level of tangles.

Let $\T(S)$ and $\ST(S)$ be the set of tangles and singular tangles, respectively, with underlying skeleton $S$. Define a linear map
\[
s:\kk[\ST(S)]\longrightarrow \kk[\T(S)]
\]
by successive application of the Vassiliev skein relation:
\[
\tik{\sing[->]{0}{0}}\longmapsto\tik{\positive[->]{0}{0}}-\tik{\negative[->]{0}{0}}
\]
It allows to identify the space of singular tangles with a subspace of $\kk[\T(S)]$. Denote by $\I_n$ the subspace of $\kk[\T(S)]$ generated by (the images of) singular tangles with at least $n$ singular points. Then the sequence
\[
\I_0=\kk[\T(S)]\supset \I_1 \supset \I_2 \dots
\]
induces a filtration on $\kk[\T(S)]$. This filtration is clearly compatible with the composition of tangles.
\begin{defi}
An invariant $f$ is a Vassiliev, or finite type, invariant of degree at most $n$ if 
\[
f(\I_{n+1}) = 0.
\]
\end{defi}
Hence, $f$ induces a well defined element of the vector space $\hom(\kk[\T(S)]/\I_{n+1},\kk)$. It follows easily from the definition that the value of a degree $n$ invariant on a tangle with $n$ singularity does not change if one replace an overcrossing by an undercrossing and \emph{vice versa}. Hence, it is sensitive only to the combinatorial information given by the positions of singularities. This information can be encoded into a chord diagram, of which we recall the definition:
\begin{defi}
A \emph{chord diagram} with underlying skeleton $S$ is a choice of $n$ pairs of points of $S\backslash \partial S$ called chord endpoints, together with an element of $\Z/2\Z$ attached to each component of $S$ called the residue. 
\end{defi}
Composition of chord diagrams is defined as for tangles. The residue on a component $k$ of the composition of two diagrams $C,C'$ is given by:
\[
\sum_i r(\alpha_i)+\sum_j r(\beta_j)+\sum_{i<i'} r(\alpha_i,\alpha_{i'})+\sum_{j<j'} r(\beta_j,\beta_{j'})
\]
where $\alpha_i,\beta_j$ are the components of $C,C'$ contained in $k$, $r(\alpha_i),r(\beta_j)$ their residue and $r(\alpha_,\alpha_{i'})=0$ if $\alpha_i,\alpha_{i'}$ can be embedded in an horizontal strip in such a way that the order of their bottom and top endpoints are preserved and without intersecting each other, and 1 otherwise.

Let $D_0(S)$ be the $\kk$ vector space generated by chord diagrams with underlying skeleton $S$. It is graded by the number of chords.
\begin{prop}
Given a chord diagram $C$ with skeleton $S$, there exists a singular tangle $T_C$ on the same skeleton such that:
\begin{itemize}
\item  two points of $S$ are joined by a chord in $C$ if and only if these are the pre-images of a singular point of $T_C$;
\item the residue on each component of $C$ matches the residue modulo 2 of the number of "twists" of the framing of the corresponding component of $T_C$. 
\end{itemize}
Every singular tangle is the realization of some diagram, and this realization is unique up to crossings flips.
\end{prop}
\begin{rmk}
	The residue reflects the fact that flipping crossings preserve the parity of the number of twists. Another way of dealing with the framing is to introduce a notion of singular framing, as in~\cite{Yetter2001}.
\end{rmk}
Hence, there exists a well defined surjective map from the space of diagrams with $n$ chords to $\I_n(S)/\I_{n+1}(S)$, and therefore a well defined degree preserving linear surjective map
\[
r:D_0(S) \longrightarrow \gr \kk[\T(S)]=\bigoplus_{n\geq 0} \I_n(S)/\I_{n+1}(S)
\]
Dualizing, we get a map $f\mapsto \tilde f$ from the space of finite type invariants to the graded dual of $D_0(S)$.
\begin{defi}
The linear map $\tilde f$ is called the \emph{symbol} of $f$.
\end{defi}
Let $A(S)$ be the quotient of $D_0(S)$ by the following relations:
\begin{align}\label{4T}\tag{4T}
\tik{\hori{0}{0}{1} \straight{2}{0}\hori[->]{0}{1}{2}\straight[->]{1}{1}}+\tik{\hori{0}{0}{1} \straight{2}{0}\hori[->]{1}{1}{1}\straight[->]{0}{1}}
=\tik{\hori[->]{0}{1}{1} \straight[->]{2}{1}\hori{0}{0}{2}\straight{1}{0}}+\tik{\hori[->]{0}{1}{1} \straight[->]{2}{1}\hori{1}{0}{1}\straight{0}{0}}
\end{align}

Let $\widehat{A(S)}$ be the degree completion of $A(S)$, and 
\[
\widehat{\kk[\T(S)]}=\lim_{\leftarrow} \kk[T(S)]/\I_{n}
\]
The fundamental theorem of Vassiliev invariants can be stated as follows:

\begin{thm}[\cite{Bar-Natan1997,Cartier1993,Kassel1998,Kontsevich1993,Le1995,Piunikhin1995}]
For each skeleton $S$, there exists a linear map
\[
Z_S: \kk[\T(S)]\longrightarrow \widehat{A(S)}
\]
compatible with the composition of tangles, defined over $\kk=\Q$ and such that if $C$ is a chord diagram and $T_C$ a realization of it, then
\[
Z_S(T_C)=C+\text{\emph{higher degree terms}}.
\]
In particular, it induces an isomorphism of filtered vector spaces.
\[
\widehat{\kk[\T(S)]}\cong \widehat{A(S)}
\]
\end{thm}
This map is called a universal Vassiliev invariant, since for any finite type invariant $f$ one has
\[
f=\tilde f\circ Z_S +\text{ finite type invariants of lower degree.}
\]

\subsection{Tangles in the solid torus}
The definition of tangles in the solid torus is the obvious analog of that of ordinary tangles. However, it will be more convenient to think again of $\{0\}\times  [0,1]$ as an additional fixed strand, hence to take the following equivalent definition:
\begin{defi}
Let $I$ be an interval. A framed oriented tangle in the solid torus (or B-tangle) with underlying skeleton $S$ is a smooth embedding of $I\sqcup S$ into $\C\times[0,1]$ such that $I$ is mapped to $\{0\}\times [0,1]$, with the same conditions and the same definition of framing as for usual tangles.
\end{defi}
\begin{rmk}
Thinking of $\C\times \{0\}$ as an additional, fixed strand makes definitions simpler, but there is no need to put an orientation or a framing on it.
\end{rmk}
It allows to define easily a generalization of the notion of singular tangle, by allowing singular crossing with the distinguished strand as well.
\begin{defi}
A singular B-tangle with skeleton $S$ is a smooth map $I\sqcup S \rightarrow \C\times [0,1]$ mapping $I$ to $\{0\}\times [0,1]$ with finitely many transverse double-points.
\end{defi}

Let $\T_B(S)$ (resp. $\ST_B(S)$) be the set of B-tangles (resp. singular B-tangles) with skeleton $S$.
\subsection{$N$-filtration}
Following Section~\ref{sec:motiv}, we introduce a way of resolving singularities of B-tangles. The main difference is that we have to take care of the framing. For reasons which will become apparent later, we will assume that the additional Artin generator of $B_n^1$ has a non-trivial framing.

Let $N\geq 1$, define a linear map
\[
s:\kk[\ST_B(S)]\longrightarrow \kk[\T_B(S)]
\]
by successive application of the usual Vassiliev skein relation for singularities between ordinary components
and the following rule for singularities with the distinguished component
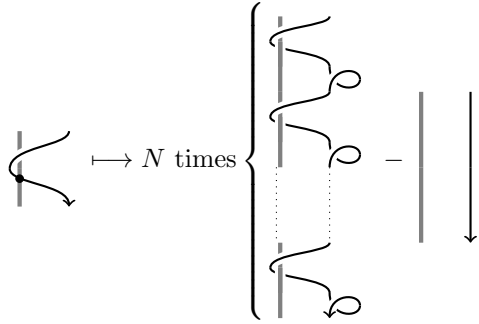
\[
\tik{\staoo[->]{0}} \ \longmapsto N \text{ times}\left\{\tik{\taotet{0}\taotet{1}\sdot{0}{2}\sdot{-0.7}{2}\taotet[->]{3}}\right.\ -\ \tik{\fixed{0}\straight{0}{0}\fixed{1}\straight[->]{0}{1}}
\]
Let $\I_{n,N}(S)$ be the image by $s$ of singular B-tangles with at least $n$ singularities of either type.
Let $f$ be an invariant, i.e. a linear map $\kk[\T_B(S)] \rightarrow \kk$. 
\begin{defi}
The map $f$ is said to be of type or degree $n$ with respect to the $N$-filtration if
\[
f(\I_{n+1,N}(S))=0\text{ and } f(\I_{n,N}(S)) \neq 0
\]
\end{defi}

\section{Examples coming from cyclotomic Hecke algebras}\label{sec:hecke}
In~\cite{Lambropoulou1999}, Lambropoulou constructed a family of Markov traces on cyclotomic Hecke algebras, of which we recall the definition: let $q, \mathbf{u}=(u_1,\dots,u_N)$ be a set of indeterminates and define $H_n(q,\mathbf u)$ as the quotient of the group algebra $\kk[q^{\pm},\mathbf u][B_n^1]$ by the relations:
\begin{align}
\sigma_i^2&=(q-1)\sigma_i+q&\prod_{i=1}^N (\tau - u_i)=&0
\end{align}
\begin{thm}[Lambropoulou]
For any choice of $z\neq 0, s_k \in \kk, k=1\dots N$, there exists $\lambda \in \kk$ and a trace
\[
tr:\bigcup_{n\geq 1} H_n(q,\mathbf u)\longrightarrow \kk[q^{\pm},\mathbf u]
\]
such that the composition of the renormalized projection
\[
\begin{array}{rcl}
B_n^1&\longrightarrow & H_n(q,\mathbf u)\\
\tau &\longmapsto & \tau\\
\sigma_i &\longmapsto & \sqrt{\lambda} \sigma_i
\end{array}
\]
with $tr$ yields a knot invariant.
\end{thm}
Now, we prove the following:
Let $\chi$ be the invariant constructed above, $\zeta$ a primitive $N$th root of unity, and assume that $\kk(\zeta)\subset \kk$. Set
\begin{align*}
q&=e^{\h} & u_i&=\zeta^i q.
\end{align*}
Hence, $\chi$ specializes to an invariant 
\[
\chi=\sum_{m\geq 0} c_m \h^m
\]
\begin{prop}
The $m$th coefficient $c_m$ of the above expansion is a $m$th order $N$-finite type invariant.
\end{prop}
\begin{proof}
The following relation holds in $H_n(q,(q\zeta,\dots,q\zeta^{N-1},q))$:
\begin{align}\label{eq:skein1}
\sigma_i-\sigma_i^{-1}&=(q-1)(\sigma_i^{-1}+1)
\end{align}
On the other hand we have
\[
0=\prod_{i=1}^N (\tau-\zeta^iq)=\tau^N-q^N
\]
therefore:
\begin{align}\label{eq:skein2}
\tau^N-1=q^N-1
\end{align}
Setting $q=e^{\h}$, the right hand sides of \eqref{eq:skein1} and \eqref{eq:skein2} are divisible by $\h$, and it clearly remains true if one compose both sides on the left and on the right by some braids $\alpha,\beta$.

It follows that the extension of $\chi$ to a singular braid with $m+1$ singularities is divisible by $\h^{m+1}$ meaning that its $m$th coefficient is 0.
\end{proof}

\section{$N$-Chord diagrams}\label{sec:ndiag}
\subsection{Definition and realization}

\begin{defi}
A $N$-\emph{chord diagram} $C$ with underlying skeleton $S$ is the data consisting of:
\begin{itemize}
\item a list $l$ of pairs of points $I\sqcup S \backslash \partial (I\sqcup S)$ called chords endpoints such than at least one of the two points of each pair is on $S$ (in other words, no chords has its two endpoints on $I$).;
\item a list of points on $S \backslash (l \cup \partial S)$ labelled by an element of $\Z/N\Z$ modulo the following rules:
\begin{align*}
\tik{\straight{0}{0} \node[point,label=left:$0$] at (0,-0.5) {};} &=\tik{\straight{0}{0}}& \tik{\straight{0}{0}\node[point,label=left:$a$] at (0,-0.3) {};\node[point,label=left:$b$] at (0,-0.7) {};}=\tik{\straight{0}{0}\node[point,label=left:$a+b$] at (0,-0.5) {};}
\end{align*}
\item for each component of $S$, an element of $\Z/2\Z$ called the residue.
\end{itemize}
\end{defi}
\begin{defi}
	Let $T$ be a B-tangle. A \emph{special point} of $T$ is either a endpoint of a component of $T$, or a singularity between two ordinary strand.
\end{defi}

Let $T$ be a B-tangle. By applying a small perturbation if necessary, we can assume that all  singularities that involve two ordinary strands belong to $(\C-(-\infty,0])\times [0,1]$. Let $\gamma \subset T$ be either a path in $T$ between two special points or a closed component of $T$ without ordinary singularities. Let $\gamma'$ be the path obtained by removing all singularities between $\gamma$ and $\{0\}\times [0,1]$ with the help of  applying the following move:

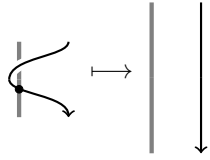
\[
\tik{\staoo[->]{0}} \ \longmapsto \ \tik{\fixed{0}\straight{0}{0}\fixed{1}\straight[->]{0}{1}}
\]
Join the ends of $\gamma'$ by a curve which does not cross $(-\infty,0]\times [0,1]$ and denote by $[\gamma]$ the element of $\pi_1(\C^*)=\Z$ obtained in this way. Since $\C-(-\infty,0]$ is simply connected it does not depend on the choice of the closure of $\gamma'$.
\begin{prop}
Let $C$ be an $N$-chord diagram with $k$ chords: 
\begin{enumerate}
\item there exists a B-tangle $T_C$ such that:
\begin{itemize}
\item The smooth map $I\sqcup S\mapsto T_C$ can be chosen in such a way that two points on $C$ are joined by a chord if and only if these are the pre-image of a given singularity;
\item if $\gamma$ is a path in $T_C$ between two usual singularities as above, then the sum of the labelling on its preimage in $C$ is equal to $[\gamma]$ modulo $N$;
\item the residue on each component is equal to the residue modulo 2 of the difference between the total number of twists of the corresponding component of $T_C$ and its representative in $\Z$.
\end{itemize}
\item The map $C\mapsto T_C$ is well defined modulo $\I_{k+1,N}(S)$.
\item Every singular tangle is $T_C$ for some diagram $C$.
\end{enumerate}
\end{prop}
\begin{proof}
(a) $T_C$ can be constructed in the obvious way: "contract" each chord in order to create singular points, add for each label the required number of loop around 0 and the same number of twists, and add a twist on each component for which the residue does not match.\\
(b) if $T_1,T_2$ are associated to $C$, one can go from $T_1$ to $T_2$ by flipping crossings and applying the following move: 
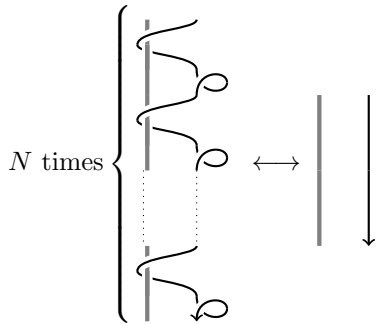
\[
N \text{ times}\left\{\tik{\taotet{0}\taotet{1}\sdot{0}{2}\sdot{-0.7}{2}\taotet[->]{3}}\right.\ \longleftrightarrow\ \tik{\fixed{0}\straight{0}{0}\fixed{1}\straight[->]{0}{1}}
\]
since the chord diagram determines the representative modulo $N$ of each path between two singular points.\\
(c) clear.
\end{proof}
Here is an example of a $N$-diagram and a realization of it. The knot is given the blackboard framing, hence its twist number is 1, and since there is one negative loop around 0, the residue of the diagram is $1-(-1)\mod 2=0$.
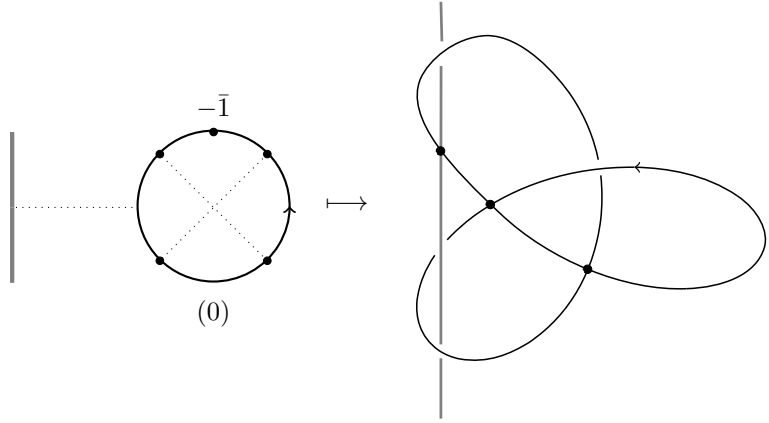
\[
\tik{\fixed{0}\fixed{-1}
\draw[dotted] (-0.6,0)--(1,0);
\begin{scope}[xshift=2cm]
\draw[->,thick] (1,0) arc (-1:360:1cm);
\node[point] at (45.0:1) {};
\node[point] at (135.0:1) {};
\draw[dotted] (45.0:1)--(225.0:1);
\node[point] at (225.0:1) {};
\node[point] at (315.0:1) {};
\node[point,label=above:$-\bar 1$] at (90.0:1) {};
\node at (-90.0:1.4) {$(0)$};
\draw[dotted] (315.0:1)--(135.0:1);
\end{scope}
}
\quad \longmapsto \quad
\begin{tikzpicture}[y=0.80pt,x=0.80pt,yscale=-1, inner sep=0pt, outer sep=0pt,baseline=(current bounding box.center)]
	\draw[normal,line join=miter,line cap=butt,miter
  limit=4.00,line width=0.577pt] (104.8056,283.2845) .. controls
  (105.6117,290.4398) and (105.5017,297.4187) .. (104.7924,304.0701) .. controls
  (101.8912,327.3802) and (91.1822,348.4825) .. (71.1477,362.1223) .. controls
  (62.6788,367.3317) and (52.2051,372.1981) .. (38.5925,369.7982) .. controls
  (35.3328,369.1344) and (31.9994,367.8136) .. (29.0846,365.8620) .. controls
  (19.6821,359.2365) and (18.2005,349.6893) .. (18.8076,342.4938) .. controls
  (19.7551,334.4676) and (22.7812,327.3248) .. (27.1440,320.9399);
  \draw[normal,line join=miter,line cap=butt,miter
  limit=4.00,line width=0.577pt] (32.9550,313.6543) .. controls
  (36.7206,309.5420) and (41.0095,305.7856) .. (45.5573,302.3404) .. controls
  (64.6582,288.4218) and (90.1965,279.1750) .. (120.4917,279.4216);
  \draw[<-,normal,line join=miter,line cap=butt,miter
  limit=4.00,line width=0.577pt] (120.4917,279.4216) .. controls
  (135.8414,279.7921) and (153.8698,282.7182) .. (168.7029,292.9211);
  \draw[normal,line join=miter,line cap=butt,miter
  limit=4.00,line width=0.577pt] (168.7029,292.9211) .. controls
  (178.5592,299.5760) and (183.6268,309.9307) .. (181.3186,317.4103) .. controls
  (177.5533,330.2362) and (161.8079,336.0911) .. (145.0538,336.7016);
  \draw[normal,line join=miter,line cap=butt,miter
  limit=4.00,line width=0.577pt] (145.0538,336.7016) .. controls
  (116.1840,337.9390) and (80.7805,323.3999) .. (57.6621,301.1606) .. controls
  (36.3154,281.0074) and (12.4672,256.3990) .. (20.5063,237.5807) .. controls
  (24.1321,229.2004) and (36.3324,217.9716) .. (51.3598,217.3597) .. controls
  (68.1007,217.1647) and (84.6982,236.2420) .. (91.3522,246.2477) .. controls
  (97.8258,256.1081) and (101.6650,266.1119) .. (103.6448,275.8893);
  \draw[zero,line join=miter,line cap=butt,miter limit=4.00,line
  width=1.024pt] (30.0000,201.4210) -- (30.4574,220.0772)(30.0757,231.6727) --
  (30.0000,282.1996) -- (30.0000,362.9783)(30.0000,369.4995) --
  (30.0000,397.9951);
\draw[cm={{1.3636389,0.0,0.0,1.3636389,(-187.86377,-51.728387)}},fill=black]
  (158.2202,237.6063)arc(-194.666:141.631:1.454);
\draw[cm={{1.3636389,0.0,0.0,1.3636389,(-119.05709,4.0391916)}},fill=black]
  (158.2202,237.6063)arc(-194.666:141.631:1.454);
\draw[cm={{1.3636389,0.0,0.0,1.3636389,(-164.5699,-26.482272)}},fill=black]
  (158.2202,237.6063)arc(-194.666:141.631:1.454);
\end{tikzpicture}
\]
Let $D_N(S)$ be the $\kk$-linear span of $N$-diagrams with skeleton $S$. Thanks to the above proposition, there is a well defined surjective map
\[
D_N(S)\longrightarrow \bigoplus_{k=0}^{\infty} \I_{k,N}(S)/\I_{k+1,N}
\]
By duality, we have the following:
\begin{cor}
Let $f$ be an invariant of $N$-type $k$. Then $f$ induces a map
\[
\tilde f:D_N(S)\longrightarrow \kk
\]
which we call its \emph{symbol}.
\end{cor}
\subsection{$N$T and labelled 4T relations}
\begin{prop}\label{prop:NT}
Let $f$ be an $N$-finite type invariant. Then its symbol $\tilde f$ descends to a function on the quotient of $D_N(S)$ by the following relations:
\begin{align}\label{eq:nat}\tag{Nat}
\tik{\hori{0}{0}{1} \node[point, label=right:$a$] at (1,0) {};\node[point, label=left:$a$] at (0,0) {};}&=\tik{\hori{0}{0}{1} \node[point, label=right:$a$] at (1,-1) {};\node[point, label=left:$a$] at (0,-1) {};} & \tik{\tzero{0}{1}\node[point, label=right:$a$] at (1,0) {};}&=\tik{\tzero{0}{1}\node[point, label=right:$a$] at (1,-1) {};}
\end{align}
\begin{align}\label{eq:lab4T}\tag{labelled 4T}
\tik{\labhori{0}{0}{1}{a}{b} \straight{2}{0}\labhori[->]{0}{1}{2}{\ }{-a-b}\straight[->]{1}{1}}+\tik{\labhori{0}{0}{1}{a}{-a} \straight{2}{0}\labhori[->]{1}{1}{1}{b}{-b}\straight[->]{0}{1}}
=\tik{\labhori[->]{0}{1}{1}{\ }{-a}\straight[->]{2}{1}\labhori{0}{0}{2}{a+b}{-b}\straight{1}{0}}+\tik{\labhori[->]{0}{1}{1}{a}{-a} \straight[->]{2}{1}\labhori{1}{0}{1}{b}{-b}\straight{0}{0}}
\end{align}

\begin{align}\tag{NT1}\label{eq:nt1}
\tik{\tzero{0}{1}\straight{2}{0}\tzero[->]{1}{2}\straight[->]{1}{1}}+\sum_a \tik{\tzero{0}{1}\straight{2}{0}\labhori[->]{1}{1}{1}{a}{-a}\fixch{1}}
&=\tik{\tzero[->]{1}{1}\straight[->]{2}{1}\tzero{0}{2}\straight{1}{0}}+\sum_a \tik{\tzero[->]{1}{1}\straight[->]{2}{1}\labhori{1}{0}{1}{a}{-a}\fixch{0}}
\end{align}
\begin{align}\tag{NT2}\label{eq:nt2}
\tik{\labhori{1}{0}{1}{a}{-a}\fixch{0}\tzero[->]{1}{2}\straight[->]{1}{1}}+\sum_b \tik{\labhori{1}{0}{1}{a}{b-a}\fixch{0}\labhori[->]{1}{1}{1}{\ }{-b}\fixch{1}}
&=\tik{\labhori[->]{1}{1}{1}{a}{-a}\fixch{1}\tzero{0}{2}\straight{1}{0}}+\sum_b \tik{\labhori{1}{0}{1}{b}{a-b}\fixch{0}\labhori[->]{1}{1}{1}{\ }{-a}\fixch{1}}
\end{align}
\end{prop}
\begin{proof}
The two first relations are consequences of the following identities:
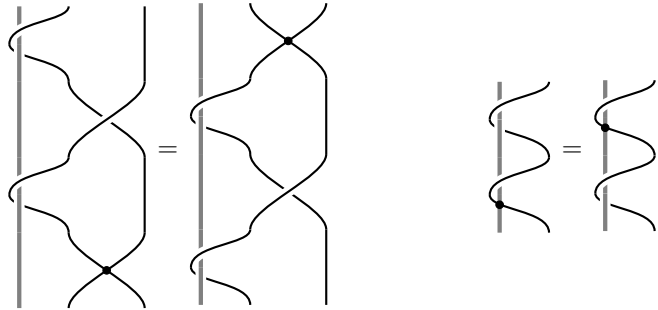
\begin{align*}
\tik{\taoo{0}\straight{1}{0}\fixed{1}\positive{0}{1}\taoo{2}\straight{1}{2}\fixed{3}\sing{0}{3}}
&=\tik{\fixed{0}\sing{0}{0}\taoo{1}\straight{1}{1}\fixed{2}\positive{0}{2}\taoo{3}\straight{1}{3}}
&
\tik{\taoo{0}\staoo{1}} &=\tik{\staoo{0}\taoo{1}}
\end{align*}
We will prove only~\eqref{eq:nt1}, the proof of the second one is similar, and let us assume that $N=3$ for simplicity. We start from the following equality:
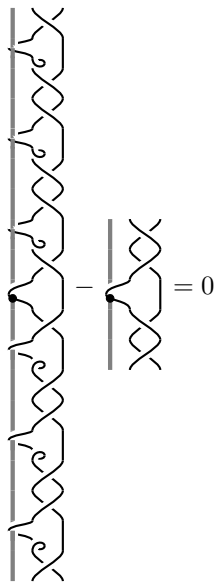
\[
\tik{
\begin{scope}[scale=0.4]
\fixed{0} \negative{0}{0}
\taotetinv{1}\straight{1}{1}
\fixed{2} \negative{0}{2}
\fixed{3} \negative{0}{3}
\taotetinv{4}\straight{1}{4}
\fixed{5} \negative{0}{5}
\fixed{6} \negative{0}{6}
\taotetinv{7}\straight{1}{7}
\fixed{8} \negative{0}{8}
\staoo{9}\straight{1}{9}
\fixed{10} \positive{0}{10}
\taotet{11}\straight{1}{11}
\fixed{12} \positive{0}{12}
\fixed{13} \positive{0}{13}
\taotet{14}\straight{1}{14}
\fixed{15} \positive{0}{15}
\fixed{16} \positive{0}{16}
\taotet{17}\straight{1}{17}
\fixed{18} \positive{0}{18}
\end{scope}
}-
\tik{
\begin{scope}[scale=0.4]
\fixed{0} \negative{0}{0}
\fixed{1} \positive{0}{1}
\staoo{2}\straight{1}{2}
\fixed{3} \negative{0}{3}
\fixed{4} \positive{0}{4}
\end{scope}
}
=0
\]
and write it as an alternating sum of a sequence of crossing flips going from one picture to the other, starting from the top:
\begin{align*}
\tik{
\begin{scope}[scale=0.4]
\fixed{0} \negative{0}{0}
\taotetinv{1}\straight{1}{1}
\fixed{2} \negative{0}{2}
\fixed{3} \negative{0}{3}
\taotetinv{4}\straight{1}{4}
\fixed{5} \negative{0}{5}
\fixed{6} \negative{0}{6}
\taotetinv{7}\straight{1}{7}
\fixed{8} \negative{0}{8}
\staoo{9}\straight{1}{9}
\fixed{10} \positive{0}{10}
\taotet{11}\straight{1}{11}
\fixed{12} \positive{0}{12}
\fixed{13} \positive{0}{13}
\taotet{14}\straight{1}{14}
\fixed{15} \positive{0}{15}
\fixed{16} \positive{0}{16}
\taotet{17}\straight{1}{17}
\fixed{18} \positive{0}{18}
\end{scope}
}\pm
\tik{
\begin{scope}[scale=0.4]
\fixed{0} \positive{0}{0}
\taotetinv{1}\straight{1}{1}
\fixed{2} \negative{0}{2}
\fixed{3} \negative{0}{3}
\taotetinv{4}\straight{1}{4}
\fixed{5} \negative{0}{5}
\fixed{6} \negative{0}{6}
\taotetinv{7}\straight{1}{7}
\fixed{8} \negative{0}{8}
\staoo{9}\straight{1}{9}
\fixed{10} \positive{0}{10}
\taotet{11}\straight{1}{11}
\fixed{12} \positive{0}{12}
\fixed{13} \positive{0}{13}
\taotet{14}\straight{1}{14}
\fixed{15} \positive{0}{15}
\fixed{16} \positive{0}{16}
\taotet{17}\straight{1}{17}
\fixed{18} \positive{0}{18}
\end{scope}
}\pm
\tik{
\begin{scope}[scale=0.4]
\fixed{0} \positive{0}{0}
\taotetinv{1}\straight{1}{1}
\fixed{2} \negative{0}{2}
\fixed{3} \positive{0}{3}
\taotetinv{4}\straight{1}{4}
\fixed{5} \negative{0}{5}
\fixed{6} \negative{0}{6}
\taotetinv{7}\straight{1}{7}
\fixed{8} \negative{0}{8}
\staoo{9}\straight{1}{9}
\fixed{10} \positive{0}{10}
\taotet{11}\straight{1}{11}
\fixed{12} \positive{0}{12}
\fixed{13} \positive{0}{13}
\taotet{14}\straight{1}{14}
\fixed{15} \positive{0}{15}
\fixed{16} \positive{0}{16}
\taotet{17}\straight{1}{17}
\fixed{18} \positive{0}{18}
\end{scope}
}\pm
\tik{
\begin{scope}[scale=0.4]
\fixed{0} \positive{0}{0}
\taotetinv{1}\straight{1}{1}
\fixed{2} \negative{0}{2}
\fixed{3} \positive{0}{3}
\taotetinv{4}\straight{1}{4}
\fixed{5} \negative{0}{5}
\fixed{6} \positive{0}{6}
\taotetinv{7}\straight{1}{7}
\fixed{8} \negative{0}{8}
\staoo{9}\straight{1}{9}
\fixed{10} \positive{0}{10}
\taotet{11}\straight{1}{11}
\fixed{12} \positive{0}{12}
\fixed{13} \positive{0}{13}
\taotet{14}\straight{1}{14}
\fixed{15} \positive{0}{15}
\fixed{16} \positive{0}{16}
\taotet{17}\straight{1}{17}
\fixed{18} \positive{0}{18}
\end{scope}
}\pm
\tik{
\begin{scope}[scale=0.4]
\fixed{0} \positive{0}{0}
\taotetinv{1}\straight{1}{1}
\taotetinv{2}\straight{1}{2}
\taotetinv{3}\straight{1}{3}
\fixed{4} \negative{0}{4}
\staoo{5}\straight{1}{5}
\fixed{6} \positive{0}{6}
\taotet{7}\straight{1}{7}
\fixed{8} \negative{0}{8}
\fixed{9} \positive{0}{9}
\taotet{10}\straight{1}{10}
\fixed{11} \positive{0}{11}
\fixed{12} \positive{0}{12}
\taotet{13}\straight{1}{13}
\fixed{14} \positive{0}{14}
\end{scope}
}\pm
\tik{
\begin{scope}[scale=0.4]
\fixed{0} \positive{0}{0}
\taotetinv{1}\straight{1}{1}
\taotetinv{2}\straight{1}{2}
\taotetinv{3}\straight{1}{3}
\fixed{4} \negative{0}{4}
\staoo{5}\straight{1}{5}
\fixed{6} \positive{0}{6}
\taotet{7}\straight{1}{7}
\fixed{8} \negative{0}{8}
\fixed{9} \positive{0}{9}
\taotet{10}\straight{1}{10}
\fixed{11} \negative{0}{11}
\fixed{12} \positive{0}{12}
\taotet{13}\straight{1}{13}
\fixed{14} \positive{0}{14}
\end{scope}
}\pm
\tik{
\begin{scope}[scale=0.4]
\fixed{0} \positive{0}{0}
\taotetinv{1}\straight{1}{1}
\taotetinv{2}\straight{1}{2}
\taotetinv{3}\straight{1}{3}
\fixed{4} \negative{0}{4}
\staoo{5}\straight{1}{5}
\fixed{6} \positive{0}{6}
\taotet{7}\straight{1}{7}
\fixed{8} \negative{0}{8}
\fixed{9} \positive{0}{9}
\taotet{10}\straight{1}{10}
\fixed{11} \negative{0}{11}
\fixed{12} \positive{0}{12}
\taotet{13}\straight{1}{13}
\fixed{14} \negative{0}{14}
\end{scope}
}\pm
\tik{
\begin{scope}[scale=0.4]
\fixed{0} \positive{0}{0}
\fixed{1}\straight{0}{1}\straight{1}{1}
\fixed{2} \negative{0}{2}
\staoo{3}\straight{1}{3}
\fixed{4} \positive{0}{4}
\taotet{5}\straight{1}{5}
\taotet{6}\straight{1}{6}
\taotet{7}\straight{1}{7}
\fixed{8} \negative{0}{8}
\end{scope}
}-
\tik{
\begin{scope}[scale=0.4]
\fixed{0} \positive{0}{0}
\fixed{1}\straight{0}{1}\straight{1}{1}
\fixed{2} \negative{0}{2}
\staoo{3}\straight{1}{3}
\fixed{4} \positive{0}{4}
\fixed{5}\straight{0}{5}\straight{1}{5}
\fixed{6} \negative{0}{6}
\end{scope}
}
=0
\end{align*}
Then, these elements pairwise leads to an additional singularity, e.g.

\[
\tik{
\begin{scope}[scale=0.4]
\fixed{0} \positive{0}{0}
\taotetinv{1}\straight{1}{1}
\fixed{2} \negative{0}{2}
\fixed{3} \negative{0}{3}
\taotetinv{4}\straight{1}{4}
\fixed{5} \negative{0}{5}
\fixed{6} \negative{0}{6}
\taotetinv{7}\straight{1}{7}
\fixed{8} \negative{0}{8}
\staoo{9}\straight{1}{9}
\fixed{10} \positive{0}{10}
\taotet{11}\straight{1}{11}
\fixed{12} \positive{0}{12}
\fixed{13} \positive{0}{13}
\taotet{14}\straight{1}{14}
\fixed{15} \positive{0}{15}
\fixed{16} \positive{0}{16}
\taotet{17}\straight{1}{17}
\fixed{18} \positive{0}{18}
\end{scope}
}-
\tik{
\begin{scope}[scale=0.4]
\fixed{0} \positive{0}{0}
\taotetinv{1}\straight{1}{1}
\fixed{2} \negative{0}{2}
\fixed{3} \positive{0}{3}
\taotetinv{4}\straight{1}{4}
\fixed{5} \negative{0}{5}
\fixed{6} \negative{0}{6}
\taotetinv{7}\straight{1}{7}
\fixed{8} \negative{0}{8}
\staoo{9}\straight{1}{9}
\fixed{10} \positive{0}{10}
\taotet{11}\straight{1}{11}
\fixed{12} \positive{0}{12}
\fixed{13} \positive{0}{13}
\taotet{14}\straight{1}{14}
\fixed{15} \positive{0}{15}
\fixed{16} \positive{0}{16}
\taotet{17}\straight{1}{17}
\fixed{18} \positive{0}{18}
\end{scope}
}
=
-\tik{
\begin{scope}[scale=0.4]
\fixed{0} \positive{0}{0}
\taotetinv{1}\straight{1}{1}
\fixed{2} \negative{0}{2}
\fixed{3} \sing{0}{3}
\taotetinv{4}\straight{1}{4}
\fixed{5} \negative{0}{5}
\fixed{6} \negative{0}{6}
\taotetinv{7}\straight{1}{7}
\fixed{8} \negative{0}{8}
\staoo{9}\straight{1}{9}
\fixed{10} \positive{0}{10}
\taotet{11}\straight{1}{11}
\fixed{12} \positive{0}{12}
\fixed{13} \positive{0}{13}
\taotet{14}\straight{1}{14}
\fixed{15} \positive{0}{15}
\fixed{16} \positive{0}{16}
\taotet{17}\straight{1}{17}
\fixed{18} \positive{0}{18}
\end{scope}
}
\]
Moreover, at each step until the last four ones, the number of loops before and after the singularity are opposite. The last four terms pairwise leads to an additional singularity with the distinguished strand, e.g.:
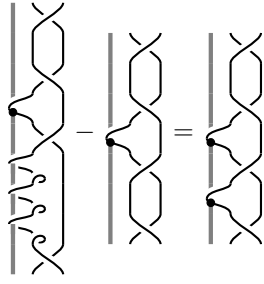
\[
\tik{
\begin{scope}[scale=0.4]
\fixed{0} \positive{0}{0}
\fixed{1}\straight{0}{1}\straight{1}{1}
\fixed{2} \negative{0}{2}
\staoo{3}\straight{1}{3}
\fixed{4} \positive{0}{4}
\taotet{5}\straight{1}{5}
\taotet{6}\straight{1}{6}
\taotet{7}\straight{1}{7}
\fixed{8} \negative{0}{8}
\end{scope}
}-
\tik{
\begin{scope}[scale=0.4]
\fixed{0} \positive{0}{0}
\fixed{1}\straight{0}{1}\straight{1}{1}
\fixed{2} \negative{0}{2}
\staoo{3}\straight{1}{3}
\fixed{4} \positive{0}{4}
\fixed{5}\straight{0}{5}\straight{1}{5}
\fixed{6} \negative{0}{6}
\end{scope}
}
=
\tik{
\begin{scope}[scale=0.4]
\fixed{0} \positive{0}{0}
\fixed{1}\straight{0}{1}\straight{1}{1}
\fixed{2} \negative{0}{2}
\staoo{3}\straight{1}{3}
\fixed{4} \positive{0}{4}
\staoo{5}\straight{1}{5}
\fixed{6} \negative{0}{6}
\end{scope}
}
\]
Hence, one get a sum of eight singular braids which are realization of the diagrams involved in \eqref{eq:nt1}.

\end{proof}

Hence, we define $B_N(S)$ as the quotient of $D_N(S)$ by the above relations.

\section{The Kontsevich--Vassiliev invariant}\label{sec:kv}
			\subsection{Ribbon categories}\label{sec:rib}
Througout this paper we only consider $\kk$-linear categories, i.e. categories whose hom spaces are $\kk$-vector spaces.  A (braided) monoidal category is the categorical analog of a (commutative) monoid, where associativity (and commutativity) holds only up to isomorphism. More precisely, a monoidal category is a category $\mC$ together with:
\begin{itemize}
 \item a bifunctor $\otimes:\mC \times \mC \longrightarrow \mC$;
 \item a natural isomorphism $\alpha:(- \ot -) \ot - \longrightarrow - \ot (- \ot -)$ (the associativity constraint);
\item a natural isomorphism $\beta: - \ot - \rightarrow - \ot^{op} -$ (the commutativity constraint, or braiding);
\item an object $I$ in $\mC$ (the unit object);
\end{itemize}
such that 
\[
\forall A \in \mC,\ I \ot A=A \ot I=A
\]

and the following diagrams commutes for any objects $A,B,C$ of $\mC$:

\begin{center}\begin{tikzpicture}[description/.style={fill=white,inner sep=2pt}]
\matrix (m) [matrix of math nodes, row sep=2em,
column sep=-0.5em, text height=1.5ex, text depth=0.25ex]
{&&(A\otimes(B\otimes C))\otimes D&& \\
((A\otimes B)\otimes C)\otimes D&&&&A\otimes((B\otimes C)\otimes D)\\
    (A\otimes B)\otimes(C\otimes D)&&&&A\otimes(B\otimes (C\otimes D))\\};
\path[->,font=\scriptsize]
(m-2-1.north) edge node[auto] {$ \alpha_{A,B,C} \ot \id_D $} (m-1-3.south west)
(m-1-3.south east) edge node[auto] {$\alpha_{A,B\ot C,D}$} (m-2-5.north)
(m-2-5.south) edge node[auto] {$ \id_A \ot \alpha_{B,C,D} $} (m-3-5.north)
(m-3-1.east) edge node[below] {$ \alpha_{A,B,C\ot D} $} (m-3-5.west)
(m-2-1.south) edge node[left] {$ \alpha_{A\ot B,C,D} $} (m-3-1.north)
;
\end{tikzpicture}\end{center}
\begin{center}\begin{tikzpicture}[description/.style={fill=white,inner sep=2pt}]
\matrix (m) [matrix of math nodes, row sep=3em,
column sep=2.5em, text height=1.5ex, text depth=0.25ex]
{& A \ot (B \ot C) & (B \ot C) \ot A &\\
(A \ot B) \ot C & && B \ot (C \ot A) \\
&(B \ot A) \ot C & B \ot (A \ot C)& \\};
\path[->,font=\scriptsize]
(m-1-2) edge node[auto] {$\beta_{A,B\ot C}$}(m-1-3)
(m-1-3) edge node[auto] {$\alpha_{B,C,A}$}(m-2-4)
(m-3-3) edge node[below right] {$\id_B\ot \beta_{A,C}$}(m-2-4)
(m-3-2) edge node[below] {$\alpha_{B,A,C}$}(m-3-3)
(m-2-1) edge node[below left] {$\beta_{A,B}\ot \id_C$}(m-3-2)
(m-2-1) edge node[auto] {$\alpha_{A,B,C}$}(m-1-2)
;
\end{tikzpicture}\end{center}
\begin{center}\begin{tikzpicture}[description/.style={fill=white,inner sep=2pt}]
\matrix (m) [matrix of math nodes, row sep=3em,
column sep=2.5em, text height=1.5ex, text depth=0.25ex]
{ & (A \ot B) \ot C  & C \ot (A \ot B) \\
A \ot (B \ot C)  & && (C \ot A )\ot B \\
&A \ot (C \ot B) & (A \ot C) \ot B 
\\};
\path[->,font=\scriptsize]
(m-1-2) edge node[auto] {$\beta_{A\ot B,C}$}(m-1-3)
(m-1-3) edge node[auto] {$\alpha_{C,A,B}^{-1}$}(m-2-4)
(m-3-3) edge node[below right] {$ \beta_{A,C}\ot \id_B$}(m-2-4)
(m-3-2) edge node[below] {$\alpha_{A,C,B}^{-1}$}(m-3-3)
(m-2-1) edge node[below left] {$\id_A \ot \beta_{B,C}$}(m-3-2)
(m-2-1) edge node[auto] {$\alpha_{A,B,C}^{-1}$}(m-1-2)
;
\end{tikzpicture}\end{center}
We also require that
\begin{align*}
\alpha_{I,A,B}&=\alpha_{A,I,B}=\alpha_{A,B,I}=\id_{A\ot B}&\beta_{I,A}=\beta_{A,I}=\id_A
\end{align*}
A strict braided monoidal category is a monoidal category for which $\alpha = \id$, and a braided monoidal category is symmetric if $\beta_{B,A} \circ \beta_{A,B}=\id_{A \ot B}$ for all $A,B \in \mC$. 

Every monoidal category is equivalent to a strict one~\cite{MacLane1963}. Let us recall from~\cite{Cartier1993} (see also~\cite[Chap.~XI.5]{Kassel1995}) the following constructive version of Mac Lane's theorem. Let $\mC^{str}$ be the category whose objects are finite sequences of objects of $\mC$. The concatenation of sequences defines an associative tensor product on $\mC^{str}$. 
For any sequence $S=(V_1,\dots,V_n)$, let $F(S)$ be its standard tensor product 
\[
F(S)=(\dots ((V_1\ot V_2)\ot V_3) \dots )\ot V_n.
\]
 Set
\[
	\hom_{\mC^{str}}(S,S'):=\hom_{\mC}(F(S),F(S'))
\]

Then $F$ induces a functor $\mC^{str}\rightarrow \mC$ which is clearly fully faithful, and since every tensor product of a finite sequence of objects in $\mC$ is canonically isomorphic to a standard one, $F$ is also essentially surjective. For any two objects of $\mC^{str}$, the associativity constraint of $\mC$ induces a canonical isomorphism
\[
F((v_1,\dots,V_k))\ot F((W_1,\dots,W_l))\cong F((V_1,\dots,V_k,W_1,\dots, W_l))
\]
which allows to define the tensor product of morphisms in $\mC^{str}$, which is therefore a strict monoidal category equivalent to $\mC$.

A strict (braided) monoidal category is called rigid if there exists a contravariant endofunctor $A \rightarrow A^*$ of $\mC$ and natural morphisms
\begin{align*}
\coev_A:I\rightarrow A \ot A^* && \ev_A:A^* \ot A \rightarrow I
\end{align*}
such that 
\[
(\id \ot \ev)\circ (\coev \ot \id) = (\ev \ot \id)\circ (\id \ot \coev) = \id
\]
If $f\in \hom(A,B)$, then $f^*\in \hom(B^*,A^*)$ is defined by
\[
f^*=(\ev_{B} \ot \id_{A^*})\circ(\id_{B^*} \ot f \ot \id_{A^*})\circ (\id_{B^*} \ot \coev_A)
\]

A ribbon structure on a strict braided monoidal category $(\mC,\id,\beta)$ is a natural automorphism $\theta$ of the identity functor satisfying:
\begin{itemize}
\item $\theta_{A \ot B}=\theta_A \ot \theta_B\beta_{B,A}\beta_{A,B}$
\item $\theta_{V^*}=(\theta_V)^*$
\end{itemize}
\subsection{The tangle category}
Let $\tang$ be the category constructed as follow. Its objects are elements of the free monoid on the set $\{+,-\}$, and the tensor product is given by concatenation. If $|w|$ is the length of a word $w$, then $\hom(w,w')$ is given by the $\kk$-linear span of isotopy classes of tangles whose skeleton $S$ satisfies $\src(S)=w,\tgt(S)=w'$. The tensor product of two morphisms $f,g$ is obtained by putting a tangle representing $g$ to the right of a tangle representing $f$ in such a way that there is no mutual linking.

\begin{rmk}
In particular, $\hom(\emptyset)$ is the algebra generated by isotopy classes of framed oriented links with composition given by the disjoint union, and $\hom(+)$ is the algebra of framed oriented links with a distinguished component with composition given by the connected sum along the distinguished component. 
\end{rmk}

A commutativity constraint $ww' \rightarrow w'w$ is given by the tangle (with the appropriate orientation) constructed by taking $|w|$ straight lines crossing a bunch of $|w'|$ straight lines:
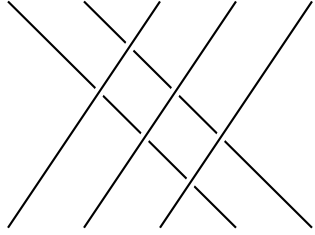
\[
\tik{
\draw[cross] (0,0) -- (3,-3);
\draw[cross] (1,0) -- (4,-3);
\draw[cross] (2,0) -- (0,-3);
\draw[cross] (3,0) -- (1,-3);
\draw[cross] (4,0) -- (2,-3);
}
\]
A ribbon structure is obtained by taking $|w|$ straight lines and applying a twist to all strands together:
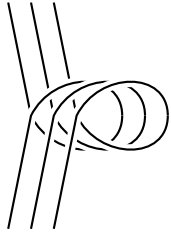
\[\tik{
\begin{scope}[scale=3]

\draw[cross] (0.5,0) -- (0.6,-0.5)..controls (0.7,-0.7) and (1,-0.7)..(1,-0.5);
\begin{scope}[xshift=0.1cm]
\draw[cross] (0.5,0) -- (0.6,-0.5)..controls (0.7,-0.7) and (1,-0.7)..(1,-0.5);
\end{scope}
\begin{scope}[xshift=0.2cm]
\draw[cross] (0.5,0) -- (0.6,-0.5)..controls (0.7,-0.7) and (1,-0.7)..(1,-0.5);
\end{scope}

\draw[cross] (1,-0.5)..controls (1,-0.3) and (0.7,-0.3).. (0.6,-0.5)--(0.5,-1);
\begin{scope}[xshift=0.1cm]
\draw[cross] (1,-0.5)..controls (1,-0.3) and (0.7,-0.3).. (0.6,-0.5)--(0.5,-1);
\end{scope}
\begin{scope}[xshift=0.2cm]
\draw[cross] (1,-0.5)..controls (1,-0.3) and (0.7,-0.3).. (0.6,-0.5)--(0.5,-1);
\end{scope}

\end{scope}
}
\]
Duality is defined by declaring that $+$ and $-$ are dual one to each other. Finally, evaluation and coevaluation are given by the following diagrams:
\begin{align*}
\tik{\up[->]{-0.5}{0}
\begin{scope}[scale=1.4]
\up[->]{-0.5}{0}
\end{scope}
\begin{scope}[scale=1.9]
\up[->]{-0.5}{0}
\end{scope}
}
&&
\tik{\ap[->]{-0.5}{-1}
\begin{scope}[scale=1.4]
\ap[->]{-0.5}{-1}
\end{scope}
\begin{scope}[scale=1.9]
\ap[->]{-0.5}{-1}
\end{scope}
}
\end{align*}
\begin{thm}[\cite{Reshetikhin1990,Shum1994}]\label{thm:RT}
The category $\tang$ is the free ribbon category on one object: for any ribbon category $\mC$ and any $V \in \mC$, there exists a structure preserving functor $F:\tang \rightarrow \mC$ mapping $+$ to $V$. 
\end{thm}

\subsection{Construction of the universal invariant}
\subsubsection{Infinitesimal braided monoidal category}
An infinitesimal braided monoidal category is a strict symmetric linear monoidal category with duality $(\mC,\id,\sigma)$ together with a natural endomorphism $t$ of the bifunctor $\ot$ satisfying:
\begin{itemize}
\item $t \circ \sigma = \sigma \circ t$
\item $t_{X,Y\ot Z}= t_{X,Y}\ot \id_Z + (\sigma_{X,Y}^{-1} \ot \id_Z) \circ (t_{X,Z} \ot \id_Y) \circ (\sigma_{X,Y}\ot \id_Z)$
\end{itemize}
Let $\cd$ be the category having the same objects as $\tang$ and whose morphisms are given by $\kk$ linear combinations of chord diagrams with composition, tensor product and duality defined as for $\tang$. Define $\widehat{\cd}$ similarly but by replacing the vector spaces of chord diagrams by their degree completion. For $w \in \cd$, let $t_w^{ij}$ be the chord diagram obtained from $\id_w$ by drawing a single chord between the $i$th and the $j$th component.

Set
\[
t_{w,w'}= \sum_{i=1}^{|w|} \sum_{j=1}^{|w'|} \epsilon_{i,j} t_{ww'}^{i,|w|+j}
\]
where $\epsilon_{i,j}=1$ if the $i$th letter of $w$ and the $j$th letter of $w'$ are identical, and $\epsilon_{i,j}=-1$ otherwise.

The braiding is defined as for $\tang$ but without distinguishing between over and undercrossings. We also declare that $t_{w,\emptyset}=t_{\emptyset,w}=\id_w$.

Let $(\mC,\sigma,t)$ be an infinitesimal braided category and define a category $\mC[[\h]]$ as follows: it has the same objects as $\mC$, and or $X,Y \in \mC$ we set
\[
\hom_{\mC[[\h]]}(X,Y)= \kk[[\h]] \ot_{\kk} \hom_{\mC}(X,Y)
\]

\begin{thm}[\cite{Cartier1993,Kassel1998}]\label{thm:univInf}
\begin{itemize}
\item $\cd$ and $\widehat{\cd}$ are infinitesimal braided categories.
\item for any infinitesimal braided category $(\mC,\sigma,t)$ and any $V \in \mC$ there exists unique functors
\[
F:\cd \longrightarrow \mC
\]
and
\[
\hat F: \widehat{\cd} \longrightarrow \mC[[\h]]
\]
such that
\[
F(+)=\hat F(+)=V
\]
and
\begin{align*}
F(t_{+,+})&=t_{V,V} & \hat F(t_{+,+})= \h t_{V,V}
\end{align*}
\end{itemize}
\end{thm}
\subsubsection{Drinfeld associators}
A \emph{horizontal} chord diagram is a chord diagram whose underlying skeleton is made of $n$ vertical intervals oriented downwards and for which each chord join two different component and does not intersect any other chord. Let $A_n$ be the space of horizontal chord diagram with $n$ components. Composition of chord diagrams turns $A_n$ into an algebra. Let
\[
t_{ij}=t^{ij}_{++\dots+}
\]
Clearly $A_n$ is generated by the $t_{ij}$'s.

Let $\hat{\mf f}_2$ be the completed free Lie algebra on two generators $a,b$ and $U(\hat{\mf f}_2)$ be its enveloping algebra. If $A$ is a complete separated filtered algebra and $x,y \in A$, then there is a unique algebra morphism sending $a$ to $x$ and $b$ to $y$. If $f \in U(\hat{\mf f}_2)$, denote by $f(x,y)$ the image of $f$ under this morphism. The algebra $U(\hat{\mf f}_2)$ has a natural structure of bialgebra with $\Delta(x)=x\ot 1+1\ot x$ and $\epsilon(x)=0$ for all $x$ in $\hat{\mf f}_2$. Recall that a \emph{group-like} element in a Hopf algebra is an element $\phi$ such that $\Delta(\phi)=\phi \ot \phi$ and $\epsilon(\phi)=1$. An element of $U(\hat{\mf f}_2)$ is group-like if and only if it is of the form $\exp(x)$ for some $x \in \hat{\mf f}_2$.
 
\begin{defi}
A Drinfeld associator, or associator for short, is a group-like element $\Phi$ of $U(\hat{\mf f}_2)$ satisfying the pentagon equation
\[
\Phi(t_{12} , t_{23} + t_{24} )\Phi(t_{13} + t_{23} , t_{34} ) = \Phi(t_{23} , t_{34} )\Phi(t_{12} + t_{13} , t_{24} + t_{34} )\Phi(t_{12} , t_{23} )
\]
in $\widehat{A}_4$ and the hexagon equation
\[
\exp(a/2)\Phi(c,a)\exp(c/2)\Phi(b,c)\exp(b/2)\Phi(a,b)=1
\]
\end{defi}

\begin{thm}[\cite{Drinfeld1990a}]
There exists an associator with rational coefficients.
\end{thm}
Let $\Phi$ be a Drinfeld associator and define natural automorphisms in $\widehat{\cd}$ by:
\begin{align*}
\alpha_{w,w',w''}&=\Phi(t_{w,w'}, t_{w',w''}) & \beta_{w,w'}&=\exp( t_{w,w'})\circ \sigma \\
\theta_w&= \exp( C_w)\circ \gamma_w& C_w&=-(\coev_w \ot d_w)\circ (t_{w,w^*}\ot \id_w) \circ (\ev_w\ot \id_w)
\end{align*}
where $\gamma_w$ is obtained from $\id_w$ by putting the residue 1 on each component.

Let $\widehat{\cd}^{str}$ be the strict monoidal category associated to $\widehat{\cd}$.

\begin{thm}[\cite{Cartier1993,Drinfeld1992,Kassel1998}]
The category $\widehat{\cd}^{str}$ is a ribbon category.
\end{thm}
Hence, there exists a functor
\[
F:\tang \longrightarrow \widehat{\cd}^{str}
\]
as in Theorem~\ref{thm:RT}.
Then the fundamental theorem of Vassiliev invariants can be stated as follows:
\begin{thm}[\cite{Bar-Natan1997,Cartier1993,Kassel1998,Kontsevich1993,Le1995,Piunikhin1995}]
The functor $F$ is a universal finite type invariant with respect to the Vassiliev filtration. Equivalently, $F$ factor through a filtration preserving isomorphism of categories
\[
\widehat{\tang} \xrightarrow{\cong} \widehat{\cd}^{str}
\]
\end{thm}
\section{Braided module categories and B-tangles invariants}\label{sec:shum}
\subsection{Braided module categories}
Let $(\mC,\alpha)$ be a monoidal category. Letting $\mC$ act on another category leads to the notion of (right) module category over a monoidal category: this is a category $\mM$ together with a bifunctor $\ot:\mM \times \mC \rightarrow \mM$ and a natural isomorphism $\gamma_{M,V,W}:(M \ot V) \ot W \longrightarrow M \ot (V \ot W)$ for $M \in \mM$ and $V,W \in \mC$ satisfying the following coherence condition:

\begin{equation}\label{eq:diagMixedPenta}
\begin{tikzpicture}[description/.style={fill=white,inner sep=2pt},scale=0.5,baseline=(current bounding box.center)]
\matrix (m) [matrix of math nodes, row sep=1.5em,
column sep=-0.5em, text height=1.5ex, text depth=0.25ex]
{&&(M\otimes(U\otimes V))\otimes W&& \\
((M\otimes U)\otimes V)\otimes W&&&&M\otimes((U\otimes V)\otimes W)\\
    (M\otimes U)\otimes(V\otimes W)&&&&M\otimes(U\otimes (V\otimes W))\\};
\path[->,font=\scriptsize]
(m-2-1.north) edge node[auto] {$ \gamma_{M,U,V} \ot \id_W $} (m-1-3.south west)
(m-1-3.south east) edge node[auto] {$\gamma_{M,U\ot V,W}$} (m-2-5.north)
(m-2-5.south) edge node[auto] {$ \id_M \ot \alpha_{U,V,W} $} (m-3-5.north)
(m-3-1.east) edge node[below] {$ \gamma_{M,U,V\ot W} $} (m-3-5.west)
(m-2-1.south) edge node[left] {$ \gamma_{M\ot U,V,W} $} (m-3-1.north)
;
\end{tikzpicture}
\end{equation}
A braided module category over a braided monoidal category $(\mC,\alpha,\beta)$ is a module category $(\mM,\gamma)$ over $(\mC,\alpha)$ together with a natural automorphism $\eta$ of the functor $\ot:\mM \times \mC \rightarrow \mM$ satisfying $\eta_{M,I}=\id_M$ for all $M\in \mM$ and such that the following diagrams commutes for $M \in \mM,\ V,W \in \mC$:

\begin{equation}\label{eq:diagOcta}
\begin{tikzpicture}[description/.style={fill=white,inner sep=2pt},scale=0.5,baseline=(current bounding box.center)]
\matrix (m) [matrix of math nodes, row sep=1.5em,
column sep=6em, text height=1.5ex, text depth=0.25ex]
{
(M \ot U)\ot V & (M \ot U) \ot V \\
M \ot (U\ot V) & M \ot (U \ot V) \\
M \ot (V\ot U) & M \ot (V \ot U) \\
(M \ot V)\ot U & (M \ot V) \ot U \\
    };
\path[->,font=\scriptsize]
(m-1-1) edge node[left] {$\gamma_{M,U,V}$} (m-2-1)
(m-2-1) edge node[left] {$\id_M \ot\beta_{U,V}$} (m-3-1)
(m-3-1) edge node[left] {$(\gamma_{M,V,U})^{-1}$} (m-4-1)
(m-4-1) edge node[below] {$\eta_{M,V} \ot \id_U$} (m-4-2)
(m-4-2) edge node[right] {$\gamma_{M,V,U}$} (m-3-2)
(m-3-2) edge node[right] {$\id_M \ot\beta_{V,U}$} (m-2-2)
(m-2-2) edge node[right] {$(\gamma_{M,U,V})^{-1}$} (m-1-2)
(m-1-1) edge node[above] {$\eta_{M\ot U,V} $} (m-1-2)
;
\end{tikzpicture}
\end{equation}

(the octagon) and
\begin{equation}\label{eq:triangle}
\begin{tikzpicture}[description/.style={fill=white,inner sep=2pt},scale=0.5,baseline=(current bounding box.center)]
\matrix (m) [matrix of math nodes, row sep=1.5em,
column sep=3em, text height=1.5ex, text depth=0.25ex]
{
(M\otimes U)\otimes V&&(M\otimes U)\otimes V&&(M\otimes U)\otimes V\\
M\otimes (U\otimes V)&&&&M\otimes (U\otimes V)\\
};
\path[->,font=\scriptsize]
(m-1-1) edge node[auto] {$ \eta_{M,U}\ot \id_V $} (m-1-3)
(m-1-3) edge node[auto] {$\eta_{M\ot U,W}$} (m-1-5)
(m-1-1) edge node[auto] {$ \gamma_{M,U,V} $} (m-2-1)
(m-2-1) edge node[below] {$ \eta_{M,U\ot V} $} (m-2-5)
(m-2-5) edge node[left] {$ \gamma_{M,U,V}^{-1} $} (m-1-5)
;
\end{tikzpicture}
\end{equation}

A module category over a strict monoidal category is strict if $\gamma=\id$. Every module category $\mM$ over a monoidal category $\mC$ is equivalent to a strict module category over $\mC^{str}$. As in Section~\ref{sec:rib}, let $\mM^{str}$ be the category whose objects are pairs of an object of $\mM$ and a finite sequence of objects of $\mC$. Concatenation of sequences defines a strict, associative module structure on $\mM^{str}$ over $\mC^{str}$. Define a functor 
\[
	F:\mM^{str}\longrightarrow\mM
\]
by
\[
F((M,(V_1,\dots,V_n))):= M\ot ((\dots ((V_1\ot V_2)\ot V_3 \dots) \ot V_n)
\]
Define
\[
	\hom_{\mM^{str}}(S,S')=\hom_{\mM}(F(S),F(S'))
\]

One shows that $\mM^{str}$ is a strict module category over $\mC^{str}$, equivalent to $(\mM,\mC)$.

\subsection{The B-tangle category}

Let $\btang$ be the category constructed as follows: its set of objects is the free right module on one generator $\bullet$ over the free monoid generated by $\{+,-\}$. Hence they are of the form $\bullet w$ where $w$ is a word on $\{+,-\}$. The morphisms set $\hom(\bullet w, \bullet w')$ is the $\kk$ vector space of formal linear combination of isotopy classes of B-tangles of type $B$ whose skeleton is of type $(w,w')$. 

Define a right action
\[
\otimes:\btang \times \tang \longrightarrow \btang
\]
by $\bullet w \otimes w' = \bullet ww'$. The tensor product of two morphisms $f \in \hom_{\btang}(a,b)$, $g \in \hom_{\tang}(c,d))$ is given again by putting $g$ on the right of $f$ in such a way that there is no mutual linking.

We have the following generalization of the Reidemeister Theorem:
\begin{thm}\label{thm:reide}
The category $\btang$ is generated by the following elementary diagrams (with all possible orientations)
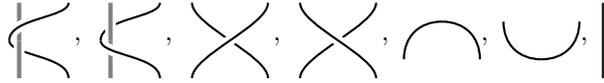
\[
\tik{\taoo{0}},\ \tik{\tainv{0}},\ \tik{\positive{0}{0}},\ \tik{\negative{0}{0}},\ \tik{\ap{0}{0}},\ \tik{\up{0}{0}},\ \tik{\straight{0}{0}}
\]
together with the following relations:
\begin{itemize}
\item Planar isotopy
\item The Turaev moves~\cite[Chap. I.3.2]{Turaev2010}
\item The reflection relation
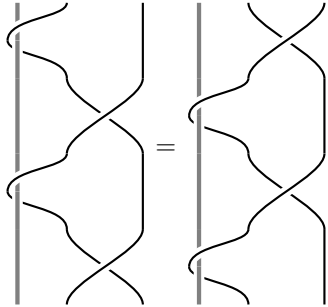
\begin{align}\label{eq:reflec}
\tik{\taoo{0}\straight{1}{0}\positive{0}{1}\fixed{1}\taoo{2}\straight{1}{2}\positive{0}{3}\fixed{3}}=\tik{\positive{0}{0}\fixed{0}\taoo{1}\straight{1}{1}\positive{0}{2}\fixed{2}\taoo{3}\straight{1}{3}}
\end{align}

\item The following relations and their upside down version:
\begin{align}\label{eq:duality}
\tik{\fixed{0}\straight{0}{0}\straight{1}{0} \fixed{1}\up{0}{1}}=\tik{\taoo{0} \straight{1}{0} \fixed{1}\positive{0}{1} \taoo{2}\straight{1}{2} \fixed{3}\negative{0}{3} \fixed{4} \up{0}{4}}=\tik{\tainv{0} \straight{1}{0} \fixed{1}\negative{0}{1} \tainv{2}\straight{1}{2} \fixed{3}\positive{0}{3} \fixed{4} \up{0}{4}}
\end{align}

\end{itemize}
\end{thm}

\begin{proof}
The proof follows from~\cite[Prop. 11]{Reinhard2001} by observing that relations (16)--(19) given there are in fact consequences of~\eqref{eq:duality} above.
\end{proof}
For any $(w_0,w) \in \btang\times \tang$, define an automorphism $\tau_{w_0,w}$ of $w_0 \ot w$ by drawing $|w_0|$ straight lines, and $|w|$ lines making a complete loop around the straight lines:
\[
\tik{
\begin{scope}[xscale=3*\taosize,yscale=4]
\draw[zero] (-1,-0)--(-1,-0.5);
\draw[normal] (-0.8,-0)--(-0.8,-0.5);
\draw[normal] (-0.6,-0)--(-0.6,-0.5);
\draw[cross] (0,-0).. controls (0,-0.2) and (-1.6,-0.2)..(-1.6,-0.5);
\begin{scope}[xshift=0.2cm]
\draw[cross] (0,-0).. controls (0,-0.2) and (-1.6,-0.2)..(-1.6,-0.5);
\end{scope}
\begin{scope}[xshift=0.4cm]
\draw[cross] (0,-0).. controls (0,-0.2) and (-1.6,-0.2)..(-1.6,-0.5);
\end{scope}

\draw[cross] (-1.6,-0.5).. controls (-1.6,-0.7) and (0,-0.7)..(0,-1);
\begin{scope}[xshift=0.2cm]
\draw[cross] (-1.6,-0.5).. controls (-1.6,-0.7) and (0,-0.7)..(0,-1);
\end{scope}
\begin{scope}[xshift=0.4cm]
\draw[cross] (-1.6,-0.5).. controls (-1.6,-0.7) and (0,-0.7)..(0,-1);
\end{scope}

\draw[zerocross] (-1,-0.5)--(-1,-2);
\draw[cross] (-0.8,-0.5)--(-0.8,-2);
\draw[cross] (-0.6,-0.5)--(-0.6,-2);

\begin{scope}[xscale=2, yshift=-1cm,xshift=-0.5cm]

\draw[cross] (0.5,0) -- (0.6,-0.5)..controls (0.7,-0.7) and (1,-0.7)..(1,-0.5);
\begin{scope}[xshift=0.1cm]
\draw[cross] (0.5,0) -- (0.6,-0.5)..controls (0.7,-0.7) and (1,-0.7)..(1,-0.5);
\end{scope}
\begin{scope}[xshift=0.2cm]
\draw[cross] (0.5,0) -- (0.6,-0.5)..controls (0.7,-0.7) and (1,-0.7)..(1,-0.5);
\end{scope}

\draw[cross] (1,-0.5)..controls (1,-0.3) and (0.7,-0.3).. (0.6,-0.5)--(0.5,-1);
\begin{scope}[xshift=0.1cm]
\draw[cross] (1,-0.5)..controls (1,-0.3) and (0.7,-0.3).. (0.6,-0.5)--(0.5,-1);
\end{scope}
\begin{scope}[xshift=0.2cm]
\draw[cross] (1,-0.5)..controls (1,-0.3) and (0.7,-0.3).. (0.6,-0.5)--(0.5,-1);
\end{scope}

\end{scope}

\end{scope}
}
\]

\begin{thm}
$(\btang,\tau)$ is a strict braided module category over $\tang$.
\end{thm}
\begin{proof}
Clear.
\end{proof}
Let $(\mC,\beta,\theta)$ be a ribbon category, $V \in \mC$ and $F$ the functor
\[
\tang \longrightarrow \mC
\] 
coming from Theorem~\ref{thm:RT}.
\begin{thm}\label{thm:myshum}
Let $(\mM,\gamma)$ be a strict braided module category over $\mC$ and $M \in \mM$. There exists a unique functor
\[
G:\btang \longrightarrow \mM
\]
extending $F$ and such that $G(\bullet)=M$ and $G(\tau_{\bullet,+})=\gamma_{M,V}$.
\end{thm}
\begin{proof}
The condition above determines $G$ uniquely. Therefore, it remains to check that it is well defined, which amounts to prove that the relations of the Theorem~\ref{thm:reide} are consequence of the axioms of braided module category. The planar isotopy reduces either to properties of the tensor product or to the naturality of $\gamma$ and $\beta$. The Turaev moves are already consequences of the axioms of ribbon categories thanks to Theorem~\ref{thm:RT}.

For $M \in \mM, X,Y \in \mC$, the reflection relation reads
\[
\gamma_{M,X}\beta_{X,Y}\gamma_{M,Y}\beta_{Y,X}=\beta_{X,Y}\gamma_{M,Y}\beta_{Y,X}\gamma_{M,X}
\]
which, using the octagon equation~\eqref{eq:diagOcta}  can be written
\[
\gamma_{M,X}\gamma_{M\ot X,Y}=\gamma_{M\ot X,Y} \gamma_{M,X}
\]
and this relation holds thanks to the naturality of $\gamma$.

The right hand side of the first equality in~\eqref{eq:duality} can be written:
\[
	\tik{\taotet{0} \straight{1}{0} \fixed{1}\positive{0}{1} \taotet{2}\straight{1}{2}
	\fixed{3}\positive{0}{3}
\twistinv{0}{4}\twistinv{1}{4}\negative{0}{4.5}\negative{0}{5.5}\fixed{4}\fixed{5}\fixed{6} \up{0}{6.5}}
\]

that is:
\begin{multline*}
	(\gamma_{M,X}\ot \id_{X^*})(\id_M\ot \beta_{X,X^*})(\gamma_{M,X^*}\ot\id_X)(\id_M \ot \beta_{X^*,X})\\(\id_M\ot \theta^{-1}_X\ot \theta^-1_{X^*})(\id_M\ot (\beta_{X,X^*}^{-1}\beta_{X^*,X}^{-1}))(\id_M\ot \ev_X)
\end{multline*}
which, using using the defining axioms of braided module categories, can be written
\[
	\gamma_{M,X\ot X^*}(\id_M\ot \theta^{-1}_{X\ot X^*})(\id_M \ot \ev_X)
\]
Then, equality~\eqref{eq:duality} holds true because of the naturality of $\gamma$ and $\theta$ and the fact that
\[
\gamma_{M,1_{\mC}}=\id_M\ot \theta_{1_{\mC}}=\id_M
\]
\end{proof}
\section{Construction of the universal invariant}\label{sec:main}

\subsection{Infinitesimal braided module categories and $N$-diagrams}
Let $(\mC,\sigma,t)$ be a infinitesimal braided module category.

\begin{defi}
An $N$-infinitesimal braided module category over $\mC$ is a strict braided module category $(\mM,\gamma)$ over $\mC$ together with a natural endomorphism $u$ of $\ot_{\mM}$ such that:
\begin{align*}
\gamma^N&=\id\\
u_{M\ot X, Y}&=\sigma_{X,Y}\circ u_{M,Y}\circ \sigma_{X,Y}+\sum_{a=0}^{N-1} \gamma_{M,X}^{-a}\circ t_{X,Y}\circ \gamma_{M,X}^a\\
u_{M,X\ot Y}&=u_{M \ot X,Y}+u_{M,X}\\
\end{align*}
\end{defi}

Let $\bcd_N$ be the category having the same objects as $\btang$ and whose morphisms are given by $N$-chord diagrams modulo the relations of Proposition~\ref{prop:NT}. We define an action of $\cd$ on $\bcd_N$ as for $\btang$.

 Define an automorphism $\tau_w$ of $w$ by labelling each  strand of $\id_w$ with $\bar 1 \in \Z/N\Z$ if it is oriented downwards, and by $-\bar 1$ otherwise. It induces an automorphism of the module structure of $\bcd_N$ by:
\[
\tau_{w_0,w}:=\tau_{w_0w}
\]

For $w \in \bcd_N$ let $t^{0i}_w$ be the $N$-diagram made from $\id_w$ by drawing a single chord between the distinguished component and the $i$th component. For $w_0 \in \bcd_N,w \in \cd$, define an endomorphism $t^0_{w_0,w}$ of $w_0 \ot w$ by:
\begin{multline}
t^0_{w_0,w}=\sum_{i=1}^{|w|} \left(w_i t^{0,|w_0|-1+i}_{w_0w}+\delta_{w_i}\sum_{\alpha}\tau_{w_0w}^{-\alpha}C_{w_0w}^{(|w_0|-1+i)}\tau_{w_0w}^{\alpha}\right)\\
+ \sum_{\substack{i,i'=1\\i\neq i'}}^{|w|} t_{ww'}^{|w_0|-1+i,|w_0|-1+i'} 
+\sum_{i=1}^{|w_0|-1} \sum_{i'=1}^{|w|}\sum_{\alpha \in \Z/N\Z} \epsilon_{i,i'}\tau_{w_0w}^{-\alpha}t_{w_0w}^{i,|w_0|-1+i'}\tau_{w_0w}^{\alpha}
\end{multline}
where $w_i$ is the $i$th letter of $w$, $\delta_{w_i}=1$ if $w_i=-$ and 0 otherwise, $C_w^{(k)}$ is $C_-$ placed on the $k$th component, $\epsilon_{i,i'}=1$ if the $i$th letter of $w_0$ ($\bullet$ having index 0) and the $i'$th letter of $w$ are identical, and $\epsilon_{i,j}=-1$ otherwise. 

\begin{thm}
The category $\bcd_N$ is an infinitesimal braided module category over $\cd$. Moreover, for any pair of an infinitesimal braided module category $(\mC,\sigma,t)$ and an $N$-infinitesimal braided module category $(\mM,\gamma,u)$ over it, the choice of objects $M,V \in \mM\times\mC$  induces a functor of braided module categories $G:\bcd_N \rightarrow \mM$ such that
\begin{align*}
G(\bullet)&=M & G(+)&=V & G(-)&=V^*\\
G(t^0_{\bullet,+})&= u_{M,V} & G(t_{+,+})&=t_{V,V}
\end{align*}
\end{thm}
\begin{proof}
The fact that these elements satisfy the axiom of infinitesimal braided module category is true by construction, but we have to show that these families of morphisms define natural morphism in $\bcd_N$.

First of all, every diagram can be written as a composition and tensor product of the identity and the following elements (with all possible residue and orientation):
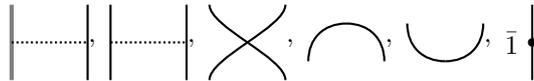
\[
\tik{\tzero{0}{1}},\ \tik{\hori{0}{0}{1}},\ \tik{\perm{0}{0}},\ \tik{\ap{0}{0}},\ \tik{\up{0}{0}} ,\tik{\zn{0}{0}{\bar 1}}
\]
hence it is enough to check the naturality on these elements.

Naturality of $\tau$ comes from~\eqref{eq:nat} which also implies that
\[
\tik{\zn{0}{0}{\bar 1}\zn[<-]{1}{0}{-\bar 1}\hori{0}{1}{1}\up{0}{2}}=\tik{\straight{0}{0}\straight[<-]{1}{0}\hori{0}{1}{1}\up{0}{2}}
\]
i.e. that $\tau_{\bullet,+-}\circ\ev_+=\ev_+$. Similar pictures shows that $\tau$ also commutes with $\coev$. That $t^0_{\bullet,++}$ commutes with the symmetric braiding is obvious. The fact that $t^0_{\bullet,++}$ and $t^0_{\bullet +,+}$ commutes with $\id_{\bullet} \ot t_{+,+}$ and $t^0_{\bullet,+}\ot \id_+$ respectively follows from~\eqref{eq:nt1} and~\eqref{eq:nt2}. That $t^0_{\bullet,+-}$ commutes with $\id_{\bullet}\ot t_{+,-}$ follows also from~\eqref{eq:nt2} since $C_-$ is central.

Finally, we have:
\begin{align*}
t^0_{\bullet,+-}\circ \ev=\tik{
\tzero[->]{0}{1} \straight[<-]{2}{0}\up{1}{1}
}
-
\tik{
\tzero[<-]{0}{2}\straight[->]{1}{0}\up{1}{1}
}
\\-
\sum_{\alpha}\tik{
\fixch{0}\straight[->]{1}{0}\straight[<-]{2}{0}\up{1}{1}
\draw[chord] (1,-0.5)--(2,-0.5);
\node[point,label=left:$-\alpha$] at (1,-0.8) {};
\node[point,label=left:$\alpha$] at (1,-0.2) {};
}
+\sum_{\alpha} \tik{ 
\fixch{0}
\begin{scope}[rotate=-90, yshift=1cm]
\draw[normal,->] (0,0) -- (1,0); \draw[chord] (0.8,0) arc (0:180:0.3); 
\node[point, label=left:$-\alpha$] at (0.5,0) {};
\node[point,label=left:$\alpha$] at (0.1,0) {};
\end{scope}
\straight[<-] {2}{0}\up{1}{1}
}
\end{align*}
which is clearly equal to 0.

The functor $G$ is uniquely determined by these requirement. Hence we have to show that it preserves the relation of Proposition~\ref{prop:NT}. The proof that the labelled 4T relation is preserved is the same as the proof of Theorem~\ref{thm:univInf} (which is Theorem~5.4 of~\cite{Kassel1998}). The first relation in~\eqref{eq:nat} follows from the naturality of $u$. The fact that the second relation of~\eqref{eq:nat} is preserved follows from the naturality of $\gamma_{M,V\ot V}$. The image of the left hand side of~\eqref{eq:nt1} is
\begin{align}\label{eq:imNT}
(u_{M,V}\ot \id_V)\circ (\sigma_{V,V} \circ u_{M,V} \circ \sigma_{V,V} +\sum_a \gamma_{M,V}^{-a}\circ t_{V,V}\circ \gamma_{M,V}^a )
\end{align}

By definition of an infinitesimal module category,
\[
\sigma_{V,V} \circ u_{M,V} \circ \sigma_{V,V} +\sum_a \gamma_{M,V}^{-a}\circ t_{V,V}\circ \gamma_{M,V}^a =u_{M\ot V,V}
\]
hence the naturality of $u$ implies that~\eqref{eq:imNT} is equal to
\[
 (\sigma_{V,V} \circ u_{M,V} \circ \sigma_{V,V} +\sum_a \gamma_{M,V}^{-a}\circ t_{V,V}\circ \gamma_{M,V}^a )\circ (u_{M,V}\ot \id_V)
\]
which is the image of the right hand side of~\eqref{eq:nt1}. The relation~\eqref{eq:nt2} is proved similarly. 
\end{proof}
Let $\mM[[\h]]$ be the category having the same objects as $\mM$, with morphisms
\[
\hom_{\mM[[\h]]}(X,Y):=\hom_{\mM}(X,Y)\ot_{\kk}\kk[[h]].
\]
Then $\h u$ and $\gamma$ induces an infinitesimal braided module structure over $\mC[[\h]]$ on $\mM[[\h]]$.
\begin{cor}
There exists a functor
\[
G:\widehat{\bcd}_N\longrightarrow \mM[[\h]]
\]
with the same properties as above, except that
\begin{align*}
G(t^0_{\bullet,+})&=\h u_{M,V}& G(t_{+,+})&=\h t_{V,V}
\end{align*}
\end{cor}
\subsection{Cyclotomic associators}
Key to the present construction is a generalization of Drinfeld associator introduced by Enriquez~\cite{Enriquez2008}.
\begin{defi}
An horizontal $N$-diagram is an $N$-chord diagram whose underlying skeleton is made of $n$ interval oriented downwards, with 0 residue and having only horizontal chords. 
\end{defi}
Let $B_{n,N}$ be the vector space of horizontal chord diagrams modulo 4T, NT1, NT2 and Nat. Again, composition of diagrams turns $B_{n,N}$ into an algebra. Let 
\[
t_{0i}=t^{0i}_{\bullet ++\dots+}
\]
and
\[
t_{ij}(a)=(\tau^{(i)})^{-a}(\id_{\bullet}\ot t_{ij})(\tau^{(i)})^a
\]
\begin{rmk}
The usual algebra of chord diagrams is the enveloping algebra of the holonomy Lie algebra $\mf t_n$ of the configuration space associated to the pure braid group~\cite{Kohno2010,Kohno1983}. Likewise, 
\[
B_{n,N}=U(\mf t_{n,N})\rtimes (\Z/N\Z)^n
\]
where $\mf t_{n,N}$ is the holonomy Lie algebra of the orbit configuration space $X_{n,N\Z}(\C^*)$ and is generated by $\{t_{0i},t_{ij}(a)|a \in \Z/N\Z\}$.
\end{rmk}
Let $\hat{\mf f}_{N+1}$ be the degree completion of the free Lie algebra on generators $\{a,b(i)|i \in \Z/N\Z\}$. If $f \in U(\hat{\mf f}_{N+1})$, $w_0 \in \bcd_N$ and $w,w' \in \cd$, set:
\[
f_{w_0,w,w'}:=f(t_{w_0,w}^0,t_{w,w'}(0),\dots,t_{w,w'}(N-1)) \in \hom_{\widehat{\bcd}_N}(w_0\ot w \ot w')
\]
\begin{defi}[\cite{Enriquez2008}]
A cyclotomic associator is a pair of group-like elements $\Psi \in U(\hat{\mf f}_{N+1})$, $\Phi \in U(\hat{\mf f}_{2})$ such that $\Phi$ is an associator and 
\[
\Psi_{\bullet,+,+}:=\Psi(t_{01},t_{12}(0),\dots,t_{12}(N-1))\in (\widehat B_{2,N})^{\times}\subset \aut_{\widehat{\bcd}_N}(\bullet++)
\] 
satisfies:
\begin{itemize}
\item the mixed pentagon equation
\[
\Psi_{\bullet,+,++}\Psi_{\bullet +,+,+}=\Phi_{+,+,+}\Psi_{\bullet,++,+}\Psi_{\bullet,+,+}
\]
in $B_{3,N}$
\item The octagon equation
\[
e^{(t_{02}+\sum_{a}t_{12}(a))/N}\tau^{(2)}=\Psi_{\bullet,+,+}^{-1} e^{t_{12}(0)/2} \Psi_{\bullet,+,+}^{0,2,1} e^{t_{02}/N}\tau^{(2)}(\Psi_{\bullet,+,+}^{0,2,1} )^{-1}e^{t_{12}(0)/2}\Psi_{\bullet,+,+}
\]
where
\begin{align*}
\Psi_{\bullet,+,+}^{0,2,1} &=\Psi(t_{02},t_{21}(0),\dots,t_{21}(N-1))\\
&=\Psi(t_{02},t_{12}(0),t_{12}(-1),\dots,t_{12}(1-N))
\end{align*}
\end{itemize}
\end{defi}
Let $\cyc(N,\kk)$ be the set of cyclotomic associators with coefficients in $\kk$.
\begin{thm}[\cite{Enriquez2008}]
For any field $\kk$ of characteristic 0, $\cyc(N,\kk)\neq \emptyset$.
\end{thm}
Let $\mf f_n^{(k)}=\mf f_n/J_k$ where $J_k=\{x \in \mf f_n|x\equiv 0 \mod \deg k\}$. The following result shows that cyclotomic associators can be constructed recursively. It is implicit in~\cite{Enriquez2008} but can be proved by reproducing the proof of~\cite[Prop. 5.8]{Drinfeld1990a}.
\begin{prop}
Let $\cyc^{(k)}(N,\kk)$ be the set of group-like elements in $U(\mf f_{N+1}^{(k)}(\kk))\times U(\mf f_{2}^{(k)}(\kk))$ satisfying the defining equations of cyclotomic associators modulo degree $k$. Then the natural map
\[
\cyc^{(k+1)}(N,\kk)\longrightarrow \cyc^{(k)}(N,\kk)
\]
is surjective.
\end{prop}
\begin{proof}[Proof]
The set $\cyc(N,\kk)$ is a torsor under the action of a certain pro-algebraic group $GRTM(N,\kk)$ which is a subset of $\exp(\hat{\mf f}_{N+1}(\kk))\times \exp(\hat{\mf f}_{2}(\kk))$. Hence one can define the algebraic group $GRTM(N,\kk)^{(k)}$ and show that its action on $\cyc(N,\kk)^{(k)}$ is free and transitive by repeating the proof of~\cite[Thm 7.10]{Enriquez2008}. Therefore, the results follows from the surjectivity of the natural map
\[
GRTM(N,\kk)^{(k+1)}\longrightarrow GRTM(N,\kk)^{(k)}
\]
which is proved in~\cite[Lemma A.3]{Enriquez2012}.
\end{proof}
\begin{rmk}
It is not known if the natural map from $\cyc(N,\kk)$ to the set of associators which send $(\Phi,\Psi)$ to $\Phi$ is surjective, i.e. if one can freely choose the underlying Drinfeld associator. It follows from~\cite[Section 11]{Enriquez2008} that this statement is implied by the generating part of the Deligne-Drinfeld conjecture on the Lie algebra $\mf{grt}$.
\end{rmk}
\subsection{The main construction}
Let $(\Psi,\Phi) \in \cyc(N,\kk)$. Since $\Phi$ is an associator, it leads to a ribbon structure on $\widehat{\cd}$ as in Section~\ref{sec:kv}. Define a natural automorphism in $\widehat{\bcd}_N$ by:
\begin{align*}
\gamma_{w_0,w'}=e^{t^0_{w_0,w}}(\id_{w_0}\ot\tau_w)
\end{align*}
\begin{thm}
The automorphisms $\Psi_{w_0,w,w'}$ and $\gamma_{w_0,w}$ turn $\widehat{\bcd}_N$ into a braided module category over the ribbon category $\widehat{\cd}$.
\end{thm}
\begin{proof}
The mixed pentagon equation and the octagon equation are, respectively, direct translations of the commutative diagrams~\eqref{eq:diagMixedPenta} and~\eqref{eq:diagOcta}. Hence it remains to prove that the diagram~\eqref{eq:triangle} commutes. Since $\bcd_N$ is an infinitesimal braided module category, 
\begin{align}\label{temp1}
t^0_{\bullet,++}=t^0_{\bullet +,+}+t^0_{\bullet,+}\ot \id_+ 
\end{align}

The naturality of $t^0$ implies that the two summand of the right hand side of the above equation commutes, hence:
\begin{align}\label{temp2}
e^{t^0_{\bullet,++}}=e^{t^0_{\bullet +,+}}\circ (e^{t^0_{\bullet,+}}\ot \id_+)
\end{align}
Since $t^0_{\bullet,+}$ obviously commutes with the right hand side of~\eqref{temp1}, it also commutes with its left hand side. On the other hand, the naturality of $t^0$ also implies that $t_{+,+}$ commutes with $t^0_{\bullet,++}$. Therefore, $\Psi_{\bullet,+,+}$ commutes with the left hand side of \eqref{temp2}, implying that
\[
\Psi_{\bullet,+,+}\circ e^{t^0_{\bullet,++}}\circ \Psi_{\bullet,+,+}^{-1}=e^{t^0_{\bullet +,+}}\circ (e^{t^0_{\bullet,+}}\ot \id_+)
\]
Finally, using relation~\eqref{eq:nat}, it follows that
\[
\Psi_{\bullet,+,+}\circ \gamma_{\bullet,++}\circ \Psi_{\bullet,+,+}^{-1}=\gamma_{\bullet +,+}\circ (\gamma_{\bullet,+}\ot \id_+)
\]
as required.
\end{proof}

Hence, by Theorem~\ref{thm:myshum} the functor 
\[
F:\tang \longrightarrow \widehat{\cd}^{str}
\]
extends to a functor
\[
G:\btang \longrightarrow \widehat{\bcd}_N^{str}
\]
The main result of this paper is the following:
\begin{thm}\label{thm:main}
The functor $G$ is a universal finite type invariant with respect to the $N$-filtration.
\end{thm}
\begin{proof}
Let $D$ be a $N$-diagram and $T$ be a singular tangle realizing $D$. We want to show that
\[
G(T)=D+\text{higher order terms}
\]
Since every tangle with $m$ singularities is the product of $m$ tangles with one singularity, it's enough to check it on these tangles. Moreover, $\Phi$ and $\Psi$ are both equal to 1 in degree 0, and since they appear only through conjugation it suffices to look at the elementary singular crossings.
Indeed,
\begin{align*}
G\left(\tik{\sing[->]{0}{0}}\right)&=\exp\left(\frac12 \tik{\hori[->]{0}{0}{1}}\right)\times \tik{\perm[->]{0}{0}}-\exp\left(-\frac12\tik{\hori[->]{0}{0}{1}}\right)\times\tik{\perm[->]{0}{0}}\\
&=\left(\tik{\perm[->]{0}{0}}-\tik{\perm[->]{0}{0}}\right)\ +\ \tik{\hori{0}{0}{1}\perm[->]{0}{1}}\ +\text{higher order terms}
\end{align*}
Likewise:
\begin{align*}
G\left(\tik{\staoo[->]{0}}\right)&=\left(\exp\left(\frac1N \tik{\tzero[->]{0}{1}}\right)\times \tik{\fixed{0}\straight[->]{0}{0}\node[point,label=left:$\bar 1$] at (0,-0.5) {};}\right)^N\ -\ \tik{\fixed{0}\straight[->]{0}{0}}\\
&=\tik{\tzero[->]{0}{1}}\ +\text{higher order terms}
\end{align*}
\end{proof}
\begin{cor}
For each skeleton $S$ and each integer $n$ there exists an isomorphism of filtered vector spaces (or algebras whenever it makes senses)
\[
\kk[\T_B(S)]/I_{n+1,N}(S) \cong B_N(S)^{\leq n}
\]
where the right hand side is the part of degree less or equal than $n$. Hence, the space of $N$-finite type invariant of degree at most $n$ is isomorphic to $(B_N(S)^{\leq n})^*$.
\end{cor}
\section{Approximating Vassiliev invariants}\label{sec:approx}
The goal of this Section is to show that finite type invariants in the usual sense are well approximated by $N$-finite type invariants. For each skeleton $S$, let $\I_n(S)$ be the $n$th piece of the Vassiliev filtration, that is the image of singular tangles having at least $n$ ordinary singularities.

\begin{prop}
We have:
\[
\bigcap_{N\geq 1} \I_{n,N}=\I_n
\]
\end{prop}
\begin{proof}
Let $\I_{n,N}^0(S)$ be the space generated by singular tangles with $n$ singularities, at least one of which is on the distinguished component. Every element of $\I_{n,N}^0(S)$ is a linear combination containing at least one tangle having at least one component whose representative in $\Z$ has an absolute value greater than $N/2$. Consequently, 
\[
\bigcap_{N\geq 0}I_{n,N}^0(S)=\{0\}
\]
\end{proof}

If $N'|N$, the projection $\Z/N\Z \rightarrow \Z/N'\Z$ induces a full functor
\[
P_{N,N'}:\bcd_N \rightarrow \bcd_{N'}
\]
by mapping $t_w^0$ to $N/N' t_w^0$ and every other elementary diagram to itself. In particular, the induced algebra map
\[
B_{2,N}\longrightarrow B_{2,N'}
\]
maps a $N$ cyclotomic associator to a $N'$ cyclotomic associator. Therefore, the collection of functors
\[
F_N:\btang \rightarrow \widehat{\bcd}_N
\]
can be chosen in such a way that all diagrams of the form
\begin{center}
\begin{tikzpicture}
\matrix (m) [matrix of math nodes, row sep=3em,
column sep=2.5em, text height=1.5ex, text depth=0.25ex]
{ \btang & & \widehat{\bcd}_N \\
& \widehat{\bcd}_{N'} & \\ };
\path[->,font=\scriptsize]
(m-1-1) edge node[auto] {$ F_N $} (m-1-3)
edge node[auto] {$ F_{N'} $} (m-2-2)
(m-1-3) edge node[auto] {$ P_{N,N'} $} (m-2-2);
\end{tikzpicture}
\end{center}
commutes.

Let
\[
	\widehat{\bcd}_{\hat\Z} = \lim_{\leftarrow} \widehat{\bcd}_N
\]
\begin{cor}
There exists a functor
\[
	\btang \longrightarrow \widehat{\bcd}_{\hat\Z}
\]
compatible with the Vassiliev filtration, whose associated graded is faithful.
\end{cor}
\section{Specializations and quantum invariants}\label{sec:quant}
\subsection{Specializations}

Let $\mf g$ be a Lie algebra over $\C$ and assume that there exists a non degenerate $t \in S^2(\mf g)^{\mf g}$. The following fact is, in a somewhat different language, due to Kontsevich~\cite{Kontsevich1993} (see also~\cite{Bar-Natan1995,Cartier1993,Kassel1998}):
\begin{prop}
The category of finite dimensional $\mf g$-module is an infinitesimal braided category, where the natural endomorphism $t_{V,W}$ is given by the action of $t$, and the symmetric braiding is given by the permutation $P:x \ot y \mapsto y \ot x$.
\end{prop}
Therefore, the choice of an associator leads to a ribbon structure on the category of $U(\mf g)[[\h]]$-module\footnote{If $A$ is an algebra, we always assume that $A[[\h]]$-modules are topologically free $\C[[\h]]$-modules of finite type.}

Let now $\sigma$ be an automorphism of $\mf g$ satisfying $\sigma^N=\id_{\mf g}$ and $(\sigma\ot\sigma)(t)=t$. The Hopf algebra structure of $U(\mf g)$ extends to $U(\mf g)\rtimes_{\sigma} \Z$ by setting
\[
\Delta(\sigma)=\sigma \ot \sigma
\]
The $\sigma$-invariance of $t$ imply that it is still a morphism in the category of $U(\mf g)\rtimes_{\sigma} \Z$-modules, which is therefore still infinitesimal braided.

 Let $\mf l=\mf g^{\sigma}$ and $\mf m=\im (\sigma-\id_{\mf g})$. Then we have a reductive decomposition
\[
\mf g=\mf l \oplus \mf m
\]
Let $t_{\mf l}$ be the image of $t$ in $S^2(\mf l)$ through the projection induced by the above decomposition. Set
\[
t_{\mf l}^{1,1}=m(t_{\mf l})
\]
where 
\[
m:U(\mf g)^{\ot 2} \rightarrow U(\mf g)
\]
is the multiplication. 
\begin{prop}
The category of finite dimensional $\mf l$-modules is a infinitesimal $N$-module category over $U(\mf g)\rtimes_{\sigma} \Z\modu$. The natural endomorphism $t^0_{M,V}$ is given by the action of $N(t_{\mf l}+\frac12 t_{\mf l}^{2,2})$ and the automorphism $\sigma_{M,V}$ by the action of $\sigma$ on $V$.
\end{prop}
\begin{proof}
The action is given by the restriction functor
\[
U(\mf g)\rtimes_{\sigma} \Z\modu\longrightarrow \mf l\modu
\] 
Since $\mf l=\mf g^{\sigma}$, $\sigma_{M,V}$ is a morphism of $\mf l$-module. Clearly, 
\[
\sum_{a\in \Z/N\Z}(\sigma \ot \id)(t)=Nt_{\mf l}
\]
Therefore,
\begin{align*}
N(\id \ot \Delta)(t_{\mf l}+\frac12 t_{\mf l}^{2,2})&=N(t_{\mf l}^{1,2}+t_{\mf l}^{1,3}+\frac12 t_{\mf l}^{2,2}+\frac12t_{\mf l}^{3,3}+t_{\mf l}^{2,3})\\
&=N(t_{\mf l}^{1,2}+\frac12 t_{\mf l}^{2,2})+N(t_{\mf l}^{1,3}+\frac12t_{\mf l}^{3,3})\\ & +\sum_{a\in \Z/NZ}(\sigma^{-a}\ot 1)t^{2,3}(\sigma^a \ot 1)
\end{align*}
as required. The second axiom is obvious.
\end{proof}
Therefore, each choice of a cyclotomic associator, a $\mf l$-module and a $U(\mf g)\rtimes_{\sigma}\Z$-module leads to a functor
\[
\btang \longrightarrow U(\mf l)[[\h]]\modu
\]
which factor through the category of $N$-chord diagrams.

\subsection{Explicit invariants}

The corresponding tangles invariants can be made explicit under some assumptions using a result of the author~\cite{Brochier2012}. Assume from now on that $\mf g$ is semi-simple and that $t$ is two times the inverse of the Killing form. Let $\mf h$ be a Cartan subalgebra of $\mf g$, and assume that $\sigma$ is given by the adjoint action of some element of the simply connected Lie group whose Lie algebra is $\mf h$, acting non-trivially on each root space of $\mf g$. 

Let $U_{\h}(\mf g)$ be the quantized enveloping algebra of $\mf g$. It is a $\C[[\h]]$-Hopf algebra generated by $(h_i,e^{\pm}_i),i=1\dots r$ where $r$ is the rank of $\mf g$, with relations which can be found in~\cite[Chap.~6.5]{Chari1994}. Let $U_{\h}(\mf h)$ be the subalgebra of $U_{\h}(\mf g)$ generated by $h_i$. 
\begin{thm}[\cite{Drinfeld1990}]
There exists an isomorphism of algebras
\[
U_{\h}(\mf g) \longrightarrow U(\mf g)[[\h]]
\]
whose restriction to $U_{\h}(\mf h)$ is the map given by $h_i \rightarrow h_i$.
\end{thm}
Therefore, there exists an equivalence of categories
\[
U(\mf g)[[\h]]\modu\longrightarrow U_{\h}(\mf g)\modu
\]
which we denote by $V[[\h]]\mapsto V_{\h}$.

\begin{thm}[\cite{Drinfeld1990,Drinfeld1990a}]
The category $U_{\h}(\mf g)\modu$ is a ribbon category which is equivalent, as a ribbon category, to $U(\mf g)[[\h]]\modu^{str}$ for any choice of associator.
\end{thm}
Therefore, for each $V \in U(\mf g)\modu$, the following diagram commutes:
\begin{center}
\begin{tikzpicture}
\matrix (m) [matrix of math nodes, row sep=3em,
column sep=2.5em, text height=1.5ex, text depth=0.25ex]
{ \tang & & U_{\h}(\mf g)\modu \\
& U(\mf g)\modu[[\h]]^{str} & \\ };
\path[->,font=\scriptsize]
(m-1-1) edge node[auto] {$ F_{V_{\h}} $} (m-1-3)
edge node[auto] {$ F_V $} (m-2-2)
(m-2-2) edge node[auto] {$ J $} (m-1-3);
\end{tikzpicture}
\end{center}

In particular, if $L$ is a framed link viewed as an endomorphism of $\emptyset$ in $\tang$, then $F_{V_{\h}}(L)=F_V(L)$.

Recall that the braided structure of $U_{\h}(\mf g)\modu$ is given by the action of some explicit invertible element $\R \in U_{\h}(\mf g)^{\ot 2}$, and its ribbon structure by some central invertible element $\theta\in U_{\h}(\mf g)$. The automorphism $\sigma$ extends to an automorphism of $U_{\h}(\mf g)$, denoted by the same symbol. The coproduct of $U_{\h}(\mf h)$ is cocommutative, meaning that $\R$ is $\mf h$-invariant which implies that $(\sigma \ot \sigma)(R)=R$. Since $\theta$ is central, the category $U_{\h}(\mf g)\rtimes_{\sigma} \Z\modu$ is still a ribbon category.

We have the following generalization of Drinfeld's result:

\begin{thm}[\cite{Brochier2012}]
Set 
\[
	E=e^{\h (t_{\mf h}+\frac12 t_{\mf h}^{2,2})}(1\ot \sigma) \in U_{\h}(\mf h)\ot U_{\h}(\mf g)\rtimes_{\sigma} \Z
\]
There exists an invertible element $\Psi \in U_{\h}(\mf h)\ot U_{\h}(\mf g)^{\ot 2}$ such that the action of $(E,\Psi)$ turns $U_{\h}(\mf h)\modu$ into a braided module category over $U_{\h}(\mf g)\rtimes_{\sigma}\Z\modu$, which is equivalent to the braided module category $U(\mf h)[[\h]]\modu$ for any choice of a cyclotomic associator.
\end{thm}

While there is no known general closed formula for $\Psi$, the corresponding invariants can still be made explicit thanks to the following
\begin{prop}[\cite{Brochier2012}]
Let $\bar \R=e^{-\h t_{\mf h}}\R$. The element $\Psi$ satisfies the following relation:
\begin{equation}\label{eq:abrr}
\Psi^{1,2,3}E^{1,2}(\Psi^{1,2,3})^{-1}=(\bar\R^{2,3})^{-1}E^{1,2}
\end{equation}
In particular, the action of 
\[
\tik{\taotet[->]{0}\straight[->]{1}{0}\straight[->]{2}{0}}\quad \dots\dots\quad \tik{\straight[->]{0}{0}}
\]
on $M_{\h}\ot (V_{\h}^{\ot n})$ is given by
\[
(\id_{M_{\h}}\ot \bar\R_{V_{\h},V_{\h}^{\ot n-1}})^{-1}E_{M_{\h},V_{\h}}\ot \id_{V_{\h}^{\ot n-1}}
\]
\end{prop}

Therefore, to the data of a $(\mf g,t,\sigma)$, a $U(\mf g)\rtimes_{\sigma}\Z$-module $V$ and a $\mf h$-module $M$, one can associate a $\eno(M[[\h]])$-valued link invariant obtained as the specialization of the universal invariant of Theorem~\ref{thm:main}, and an explicit $\eno(M_{\h})$-valued link invariant obtained from the braided module structure of $U_{\h}(\mf h)\modu$.

\begin{cor}
For any framed oriented link $L$ in the solid torus, these two invariants are equivalent. In particular, if $M$ is one dimensional, then they are actually equal.
\end{cor}

\subsection{Rational form and specializations}
Let $U_q(\mf g)$ be the rational form of $U_{\h}(\mf g)$, which is defined over the field $\Q(q)$ (see e.g.~\cite[Chap.~9.1]{Chari1994}). Let $d$ be the smallest integer such that for all $\lambda,\mu$ in the weight lattice of $\mf g$, $(\lambda,\mu)\in \frac1d \Z$ (note that $d$ is a divisor of the determinant of the Cartan matrix of $\mf g$).

\begin{thm}[\cite{Lusztig2010,Le2000,Turaev2010}]
The braiding and the ribbon element of $U_{\h}(\mf g)$ act on any module by a matrix whose coefficient are Laurent polynomial in $e^{\h/d}$. Therefore, the category of finite dimensional $U_q(\mf g)$-modules is a ribbon category, after extending the scalar from $\Q(q)$ to $\Q(\nu)$ where $\nu^d=q$. In particular, for any $U_q(\mf g)$ module $V_q$, the associated framed oriented link invariant takes values in $\Q[\nu,\nu^{-1}]$.
\end{thm}

The key fact here is that coefficients of $\bar\R$ are actual Laurent polynomial in $q$, and that $q^{t_{\mf h}}$ acts on any weight subspace $V_{\lambda}\ot V_{\mu}$ by multiplication by $q^{2(\lambda,\mu)}$.

Likewise, it was shown by the author that $\Psi$ lives in the graded completion of the localization of $U_q(\mf h)\ot U_q(\mf g)^{\ot 2}$ by some explicitly described multiplicative subset of $U_q(\mf h)^{\ot 3}$. Since $q^{t_{\h}^{1,1}/2}$ acts on any weight subspace $V_{\lambda}$ by $q^{(\lambda,\lambda)}$, we have:

\begin{prop}[\cite{Brochier2012}]
The category of finite dimensional $U_q(\mf h)$-modules is a braided module category over $U_q(\mf g)\modu$.
\end{prop}

Let $\lambda \in \mf h^*$ be an integral weight. It extends to an algebra morphism
\[
\lambda:U_q(\mf h)\longrightarrow \C(q)
\]
that is to a $U_q(\mf h)$-module structure on $\C(q)$. Therefore, for any $U_q(\mf g)$ module $V_q$, the resulting invariant for link in the solid torus belongs to $\C(\nu)$.
Finally, letting $\zeta$ be a primitive $N$th root of unity and using again relation~\eqref{eq:abrr}:
\begin{cor}
The above invariant belongs to $\Q(\zeta)[\nu,\nu^{-1}]$.
\end{cor}
\subsection{Example}
Let $\mf g=\mf{sl}_2$, and $V_{\h}=\C^2$ the fundamental $U_{\h}(\mf g)$-module. Let $\zeta$ be a $N$th root of unity and define an automorphism of $\mf g$ by setting:
\begin{align*}
	\sigma(e^{\pm})&=\zeta^{\pm 1}e^{\pm}& \sigma(h)&=h
\end{align*}
where $e^{\pm},h$ are the standard Chevalley generators of $\mf{sl}_2$. Then $V_{\h}$ can be turned into a $U_{\h}(\mf{sl}_2)\rtimes_{\sigma} \Z$-module by setting
\[
\sigma \longmapsto
\begin{pmatrix}
	\zeta^{1/2} & 0\\ 0 & \zeta^{-1/2}
\end{pmatrix}
\]
Let $\lambda \in \Z$ and $\C_{\lambda}$ be the one dimensional $U_{\h}(\mf h)$ module defined by $h\mapsto \lambda$. Then the invariant attached to the following knot with blackboard framing is 
\[
	q^{5+\lambda/2}\zeta^{1/2} + q^{1+\lambda/2}\zeta^{1/2} -  q^{-1+\lambda/2}\zeta^{1/2} - q^{-2+\lambda/2}\zeta^{-1/2}+ 2q^{\lambda/2}\zeta^{-1/2} 
\]
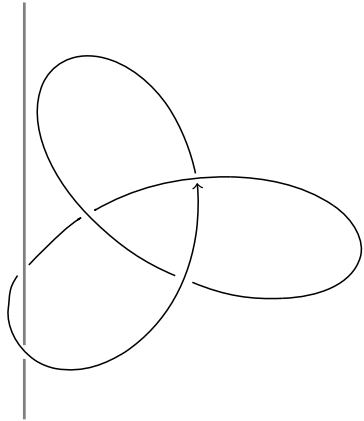
\begin{center}
\begin{tikzpicture}[y=0.80pt,x=0.80pt,yscale=-1, inner sep=0pt, outer sep=0pt]
\begin{scope}[draw,line join=miter,line cap=butt,line width=0.800pt]
\end{scope}
\draw[<-,normal,line join=miter,line cap=butt,miter
  limit=4.00,line width=0.577pt] (104.9954,282.3358) .. controls
  (105.8015,289.4911) and (105.5017,297.4187) .. (104.7924,304.0701) .. controls
  (101.8912,327.3802) and (91.1822,348.4825) .. (71.1477,362.1223) .. controls
  (62.6788,367.3317) and (52.2051,372.1981) .. (38.5925,369.7982) .. controls
  (35.3328,369.1344) and (31.9994,367.8136) .. (29.0846,365.8620) .. controls
  (19.6821,359.2365) and (15.9473,349.3954) .. (16.5545,342.1999) .. controls
  (17.5020,334.1737) and (16.6091,332.3609) .. (20.9720,325.9761);
\draw[normal,line join=miter,line cap=butt,miter
  limit=4.00,line width=0.577pt] (26.2387,320.9781) .. controls
  (30.0043,316.8658) and (41.0095,305.7856) .. (45.5573,302.3404) .. controls
  (47.7008,300.7785) and (48.2505,300.1128) ..
  (50.5536,298.6754)(50.7630,298.7452) .. controls (48.6875,299.8594) and
  (49.3169,299.4618) .. (49.5458,299.3212)(57.0631,295.1546) .. controls
  (57.7639,294.7241) and (57.0430,295.0939) ..
  (57.7581,294.6757)(56.9207,295.3037) .. controls (74.6896,284.9117) and
  (94.9592,279.2138) .. (120.4917,279.4216);
\draw[normal,line join=miter,line cap=butt,miter
  limit=4.00,line width=0.577pt] (120.4917,279.4216) .. controls
  (135.8414,279.7921) and (153.8698,282.7182) .. (168.7029,292.9211);
\draw[normal,line join=miter,line cap=butt,miter
  limit=4.00,line width=0.577pt] (168.7029,292.9211) .. controls
  (178.5592,299.5760) and (183.6268,309.9307) .. (181.3186,317.4103) .. controls
  (177.5533,330.2362) and (161.8079,336.0911) .. (145.0538,336.7016);
\draw[normal,line join=miter,line cap=butt,miter
  limit=4.00,line width=0.577pt] (145.0538,336.7016) .. controls
  (130.2717,337.3352) and (118.4233,335.7011) ..
  (102.7415,329.1465)(94.7433,326.0221) .. controls (79.7982,319.7754) and
  (68.9432,312.0128) .. (57.6621,301.1606) .. controls (36.3154,281.0074) and
  (24.6957,254.4364) .. (32.7348,235.6181) .. controls (36.3606,227.2378) and
  (45.7265,220.3201) .. (59.4729,222.9108) .. controls (73.0113,225.4913) and
  (84.6982,236.2420) .. (91.3522,246.2477) .. controls (97.8258,256.1081) and
  (102.2342,268.0094) .. (104.2141,277.7869);
\draw[zero,line join=miter,line cap=butt,miter limit=4.00,line
  width=1.024pt] (24.1290,197.1901) -- (24.1290,358.7474)(24.1290,365.2687) --
  (24.1290,393.7643);

\end{tikzpicture}

\end{center}



\begin{thebibliography}{CDM12}

\bibitem[AMR]{Andersen1998}
\textsc{J.~E. Andersen, J.~Mattes, N.~Reshetikhin}.
\newblock Quantization of the algebra of chord diagrams.
\newblock \emph{Math. Proc. Cambridge Philos. Soc.}  (\textbf{1998}).
\newblock 124(3):451--467.

\bibitem[Be]{Bezrukavnikov1994}
\textsc{R.~Bezrukavnikov}.
\newblock Koszul {DG}-algebras arising from configuration spaces.
\newblock \emph{Geom. Funct. Anal.}  (\textbf{1994}).
\newblock 4(2):119--135.

\bibitem[BF]{Bellingeri2004157}
\textsc{P.~Bellingeri, L.~Funar}.
\newblock Braids on surfaces and finite type invariants.
\newblock \emph{Comptes Rendus Mathematique}  (\textbf{2004}).
\newblock 338(2):157 -- 162.

\bibitem[BN1]{Bar-Natan1995}
\textsc{D.~Bar-Natan}.
\newblock On the {V}assiliev knot invariants.
\newblock \emph{Topology}  (\textbf{1995}).
\newblock 34(2):423--472.

\bibitem[BN2]{Bar-Natan1997}
\textsc{D.~Bar-Natan}.
\newblock Non-associative tangles.
\newblock In \emph{Geometric topology ({A}thens, {GA}, 1993)}, vol.~2 of
  \emph{AMS/IP Stud. Adv. Math.}, pp. 139--183 (Amer. Math. Soc., Providence,
  RI, \textbf{1997}).

\bibitem[Bri]{Brieskorn1971}
\textsc{E.~Brieskorn}.
\newblock Die {F}undamentalgruppe des {R}aumes der regul\"aren {O}rbits einer
  endlichen komplexen {S}piegelungsgruppe.
\newblock \emph{Invent. Math.}  (\textbf{1971}).
\newblock 12:57--61.

\bibitem[Bro]{Brochier2012}
\textsc{A.~Brochier}.
\newblock A {K}ohno–{D}rinfeld theorem for the monodromy of cyclotomic {KZ}
  connections.
\newblock \emph{Communications in Mathematical Physics}  (\textbf{2012}).
\newblock 311:55--96.
\newblock 10.1007/s00220-012-1424-0.

\bibitem[Ca]{Cartier1993}
\textsc{P.~Cartier}.
\newblock Construction combinatoire des invariants de {V}assiliev-{K}ontsevich
  des n\oe uds.
\newblock In \emph{R.{C}.{P}.\ 25, {V}ol.\ 45 ({F}rench) ({S}trasbourg,
  1992--1993)}, vol. 1993/42 of \emph{Pr\'epubl. Inst. Rech. Math. Av.}, pp.
  1--10 (Univ. Louis Pasteur, Strasbourg, \textbf{1993}).

\bibitem[CDM]{Chmutov2012}
\textsc{S.~Chmutov, S.~Duzhin, J.~Mostovoy}.
\newblock \emph{Introduction to {V}assiliev knot invariants} (Cambridge
  University Press, Cambridge, \textbf{2012}).

\bibitem[CEE]{Calaque2009}
\textsc{D.~Calaque, B.~Enriquez, P.~Etingof}.
\newblock Universal {KZB} equations: The elliptic case.
\newblock In \textsc{Y.~Tschinkel, Y.~Zarhin}, editors, \emph{Algebra,
  Arithmetic, and Geometry}, vol. 269 of \emph{Progress in Mathematics}, pp.
  165--266 (Birkhäuser Boston, \textbf{2009}).

\bibitem[CP]{Chari1994}
\textsc{V.~Chari, A.~Pressley}.
\newblock \emph{A guide to quantum groups} (Cambridge University Press,
  Cambridge, \textbf{1994}).

\bibitem[Dr1]{Drinfeld1990a}
\textsc{V.~G. Drinfeld}.
\newblock On quasitriangular quasi-{H}opf algebras and on a group that is
  closely connected with {${\rm Gal}(\overline{\bf Q}/{\bf Q})$}.
\newblock \emph{Leningrad Math. J.}  (\textbf{1990}).
\newblock 2(4):829--860.

\bibitem[Dr2]{Drinfeld1990}
\textsc{V.~G. Drinfeld}.
\newblock Quasi-{H}opf algebras.
\newblock \emph{Leningrad Math. J.}  (\textbf{1990}).
\newblock 1(6):1419--1457.

\bibitem[Dr3]{Drinfeld1992}
\textsc{V.~G. Drinfeld}.
\newblock On the structure of quasitriangular quasi-{H}opf algebras.
\newblock \emph{Funct. Anal. Appl.}  (\textbf{1992}).
\newblock 26(1):63--65.

\bibitem[EF]{Enriquez2012}
\textsc{B.~Enriquez, H.~Furusho}.
\newblock Mixed pentagon, octagon, and broadhurst duality equations.
\newblock \emph{Journal of Pure and Applied Algebra}  (\textbf{2012}).
\newblock 216(4):982 -- 995.

\bibitem[En1]{Enriquez2008}
\textsc{B.~Enriquez}.
\newblock Quasi-reflection algebras and cyclotomic associators.
\newblock \emph{Selecta Mathematica, New Series}  (\textbf{2008}).
\newblock 13:391--463.
\newblock 10.1007/s00029-007-0048-2.

\bibitem[En2]{Enriquez2005b}
\textsc{B.~Enriquez}.
\newblock Flat connections on configuration spaces and formality of braid
  groups of surfaces.
\newblock \emph{ArXiv e-prints}  (\textbf{2011}).
\newblock ArXiv.org:1112.0864.

\bibitem[Go]{Goryunov1997}
\textsc{V.~Goryunov}.
\newblock Finite order invariants of framed knots in a solid torus and in
  {A}rnold's {$J^+$}-theory of plane curves.
\newblock In \emph{Geometry and physics ({A}arhus, 1995)}, vol. 184 of
  \emph{Lecture Notes in Pure and Appl. Math.}, pp. 549--556 (Dekker, New York,
  \textbf{1997}).

\bibitem[Ha]{Habiro2000}
\textsc{K.~Habiro}.
\newblock Claspers and finite type invariants of links.
\newblock \emph{Geom. Topol.}  (\textbf{2000}).
\newblock 4:1--83 (electronic).

\bibitem[Hu]{Humbert2012}
\textsc{P.~Humbert}.
\newblock \emph{Intégrale de Kontsevich elliptique et enchevêtrements en
  genre supérieur}.
\newblock Ph.D. thesis, IRMA (\textbf{2012}).

\bibitem[Ka]{Kassel1995}
\textsc{C.~Kassel}.
\newblock \emph{Quantum groups}, vol. 155 of \emph{Graduate Texts in
  Mathematics} (Springer-Verlag, New York, \textbf{1995}).

\bibitem[Koh1]{Kohno1983}
\textsc{T.~Kohno}.
\newblock On the holonomy {L}ie algebra and the nilpotent completion of the
  fundamental group of the complement of hypersurfaces.
\newblock \emph{Nagoya Math. J.}  (\textbf{1983}).
\newblock 92:21--37.

\bibitem[Koh2]{Kohno2010}
\textsc{T.~Kohno}.
\newblock Bar complex, configuration spaces and finite type invariants for
  braids.
\newblock \emph{Topology and its Applications}  (\textbf{2010}).
\newblock 157(1):2 -- 9.

\bibitem[Kon]{Kontsevich1993}
\textsc{M.~Kontsevich}.
\newblock Vassiliev's knot invariants.
\newblock In \emph{I. {M}. {G}elfand {S}eminar}, vol.~16 of \emph{Adv. Soviet
  Math.}, pp. 137--150 (Amer. Math. Soc., Providence, RI, \textbf{1993}).

\bibitem[KT]{Kassel1998}
\textsc{C.~Kassel, V.~Turaev}.
\newblock Chord diagram invariants of tangles and graphs.
\newblock \emph{Duke Math. J.}  (\textbf{1998}).
\newblock 92(3):497--552.

\bibitem[La]{Lambropoulou1999}
\textsc{S.~Lambropoulou}.
\newblock {Knot theory related to generalized and cyclotomic {H}ecke algebras
  of type B.}
\newblock \emph{J. Knot Theory Ramifications}  (\textbf{1999}).
\newblock 8(5):621--658.

\bibitem[Le]{Le2000}
\textsc{T.~T.~Q. Le}.
\newblock Integrality and symmetry of quantum link invariants.
\newblock \emph{Duke Math. J.}  (\textbf{2000}).
\newblock 102(2):273--306.

\bibitem[Li]{Lieberum2004}
\textsc{J.~Lieberum}.
\newblock Universal {V}assiliev invariants of links in coverings of
  3-manifolds.
\newblock \emph{J. Knot Theory Ramifications}  (\textbf{2004}).
\newblock 13(4):515--555.

\bibitem[LM]{Le1995}
\textsc{T.~Q.~T. Le, J.~Murakami}.
\newblock Representation of the category of tangles by {K}ontsevich's iterated
  integral.
\newblock \emph{Comm. Math. Phys.}  (\textbf{1995}).
\newblock 168(3):535--562.

\bibitem[Lu]{Lusztig2010}
\textsc{G.~Lusztig}.
\newblock \emph{Introduction to quantum groups}.
\newblock Modern Birkh\"auser Classics (Birkh\"auser/Springer, New York,
  \textbf{2010}).
\newblock Reprint of the 1994 edition.

\bibitem[ML]{MacLane1963}
\textsc{S.~Mac~Lane}.
\newblock Natural associativity and commutativity.
\newblock \emph{Rice Univ. Studies}  (\textbf{1963}).
\newblock 49(4):28--46.

\bibitem[Pi]{Piunikhin1995}
\textsc{S.~Piunikhin}.
\newblock Combinatorial expression for universal {V}assiliev link invariant.
\newblock \emph{Comm. Math. Phys.}  (\textbf{1995}).
\newblock 168(1):1--22.

\bibitem[RH]{Reinhard2001}
\textsc{Reinhard, Häring-Oldenburg}.
\newblock Actions of tensor categories, cylinder braids and their kauffman
  polynomial.
\newblock \emph{Topology and its Applications}  (\textbf{2001}).
\newblock 112(3):297 -- 314.

\bibitem[RT]{Reshetikhin1990}
\textsc{N.~Y. Reshetikhin, V.~G. Turaev}.
\newblock Ribbon graphs and their invariants derived from quantum groups.
\newblock \emph{Comm. Math. Phys.}  (\textbf{1990}).
\newblock 127(1):1--26.

\bibitem[Sh]{Shum1994}
\textsc{M.~C. Shum}.
\newblock Tortile tensor categories.
\newblock \emph{J. Pure Appl. Algebra}  (\textbf{1994}).
\newblock 93(1):57--110.

\bibitem[Tu]{Turaev2010}
\textsc{V.~G. Turaev}.
\newblock \emph{Quantum invariants of knots and 3-manifolds}, vol.~18 of
  \emph{de Gruyter Studies in Mathematics} (Walter de Gruyter \& Co., Berlin,
  \textbf{2010}), revised ed.

\bibitem[Va]{Vassiliev1990}
\textsc{V.~A. Vassiliev}.
\newblock Cohomology of knot spaces.
\newblock In \emph{Theory of singularities and its applications}, vol.~1 of
  \emph{Adv. Soviet Math.}, pp. 23--69 (Amer. Math. Soc., Providence, RI,
  \textbf{1990}).

\bibitem[Ye]{Yetter2001}
\textsc{D.~N. Yetter}.
\newblock \emph{Functorial knot theory}, vol.~26 of \emph{Series on Knots and
  Everything} (World Scientific Publishing Co. Inc., River Edge, NJ,
  \textbf{2001}).
\newblock Categories of tangles, coherence, categorical deformations, and
  topological invariants.

\end{thebibliography}
\end{document}